\documentclass[11pt]{amsart}
\usepackage{fullpage}
\usepackage{latexsym}
\usepackage{a4}
\usepackage{amssymb}
\usepackage{amsbsy}
\usepackage{bm}
\usepackage{algpseudocode}
\usepackage{enumitem}
\newtheorem{thm}{Theorem}[section]
\newtheorem{pro}[thm]{Proposition}
\newtheorem{lem}[thm]{Lemma}

\newtheorem{cor}[thm]{Corollary}
\theoremstyle{definition}

\newtheorem{de}[thm]{Definition}
\numberwithin{equation}{section}

\newcommand{\ab}[1]{{\mathbf{#1}}}
\newcommand{\ob}[1]{{\mathbb{#1}}}

\newcommand{\N}{\Bbb{ N}}
\newcommand{\Z}{\Bbb{ Z}}

\newcommand{\ul}[1]{\underline{#1}}

\newcommand{\Con}{\mathrm{Con}}

\newcommand{\meet}{\wedge}
\newcommand{\join}{\vee}
\newcommand{\bigjoin}{\bigvee}
\newcommand{\bigmeet}{\bigwedge}

\newcommand{\HC}[1]{$(\mathrm{HC}#1)$}

\newcommand{\lcover}{\prec}

\newcommand{\algop}[2]{( {#1}, {#2} )}

\newcommand{\tup}[3]{(#1_{#2},\dots,#1_{#3})}
\DeclareMathAlphabet{\mathbfsl}{OT1}{cmr}{bx}{it}
\newcommand{\tupBold}[1]{\mathbfsl{#1}}
\newcommand{\vb}[1]{\tupBold{#1}}
\newcommand{\gb}[1]{\boldsymbol#1}

\renewcommand{\emptyset}{\varnothing}

\newcommand{\NW}[1]{\N_0^{#1}}
\newcommand{\NWC}[1]{(\N_0 \cup \{ \infty \})^{#1}}
\newcommand{\NWO}[1]{\NW{#1}}

\newcommand{\Min}{\operatorname{Min}}
\newcommand{\Max}{\operatorname{Max}}

\newcommand{\one}[2]{\mathbf{s}_{(#1,#2)}}
\newcommand{\ext}[1]{\widehat{#1}}
\title{Finite representation of commutator sequences}
\author{Erhard Aichinger}

\address{Erhard Aichinger,
Institut f\"ur Algebra,
Johannes Kepler Universit\"at Linz,
4040 Linz,
Austria}
\email{\tt erhard@algebra.uni-linz.ac.at}

\author{Neboj\v{s}a Mudrinski}
\address{Neboj\v{s}a Mudrinski,
Department of Mathematics and Informatics,
Faculty of Sciences,
University of Novi Sad,
21000 Novi Sad,
Serbia}
\email{\tt nmudrinski@dmi.uns.ac.rs}

\subjclass[2020]{08A05 (06A07, 06B99, 08A40)}
\thanks{Supported by the Austrian Science Fund (FWF):P33878 and the
Ministry of Science, Education and Technological Development of the Republic of Serbia (Grant No. 451-03-68/2022-
14/200125)}
\keywords{Lattices, sequences of commutator operations,
          higher commutators, antitone functions}
\date{\today}

\date{\today}

\begin{document}
\bibliographystyle{amsplain}
\begin{abstract}
   Several structural properties of a universal algebra
   can be seen from the higher commutators of its congruences.
   Even on a finite algebra, the sequence of higher commutator
   operations is an infinite object. In the present paper,
   we exhibit finite representations of this sequence.
\end{abstract}

\maketitle

\section{Introduction} \label{sec:intro}

In the 1970s, commutators were defined in a universal algebraic
setting \cite{Sm:MV, HH:ACIM} as a generalization of
ideal products in ring theory and the commutator subgroup in
group theory. For an arbitrary universal algebra, we rather
work with its congruence relations (these are equivalence relations that
are invariant under the operations of the algebra) than with
particular substructures, such as ideals or normal subgroups.
The \emph{commutator} of two congruences of
an algebra is another congruence of this algebra; for a
definition, see, e.g., \cite[Definition~4.150]{MMT:ALVV}. For algebras in
congruence modular varieties, this commutator operation encodes relevant
information, e.g., whether a finite algebra generates a residually
finite variety \cite[Theorem~10.15]{FM:CTFC} or whether
the logarithm of the size of the free algebra over $n$ elements
in a finitely generated variety is
bounded by a polynomial in $n$ \cite{BB:FSON, Ke:CMVW}.
In clone theory, commutators have been used to describe
local linearity properties of polynomials (cf. \cite[Lemma~3.1]{Ai:BTFS}),
to describe polynomial completeness properties
(cf. \cite{HH:ALEC, IS:PRA, KP:PCIA, Ai:OHAH, AM:TOPC}), and to distinguish
algebras with different clones of polynomial functions
\cite{AM:PCOG}. However, there is the limitation
that even in expanded groups, the polynomial functions
are often not determined by the congruence lattice and the
commutator alone. Therefore, A.\ Bulatov \cite{Bu:OTNO}
introduced higher commutators
as a finer tool to distinguish between polynomially inequivalent algebras.
His higher commutator operation associates a congruence
$[\alpha_1, \ldots, \alpha_n]$, the \emph{higher commutator} of
$\alpha_1, \ldots, \alpha_n$, with every finite sequence of
congruences of an algebra. For an algebra with congruence lattice $\ob{L}$,
these higher commutator operations form a sequence
$(f_n)_{n \in \N}$, where, for every $n \in \N$, $f_n$ is a function
from $\ob{L}^n$ to $\ob{L}$ defined by $f_n (\alpha_1, \ldots, \alpha_n)
:= [\alpha_1, \ldots, \alpha_n]$. Bulatov also states some pivotal
properties of the higher commutator operations, such as monotony
or symmetry \cite[Proposition~1]{Bu:OTNO}.
In the present note, we show that often such a sequence can be finitely
represented.

\begin{de}[cf. {\cite[p.861]{AM:SOCO}}]
  Let $\ob{L}$ be a lattice. Then an \emph{operation sequence} on $\ob{L}$
  is a sequence $(f_n)_{n \in \N}$, where, for each $n \in \N$,
  $f_n$ is a function from $\ob{L}^n$ to~$\ob{L}$.
\end{de}
If $\ob{L}$ has at least two elements, then there are at least
$2^{\aleph_0}$ operation sequences on $\ob{L}$. However, not every
operation sequence is the sequence of commutator operations of some algebra.
We will therefore restrict our attention to those operation sequences that satisfy
some additional properties that we list in the following definition.
\begin{de}[cf. {\cite{Bu:OTNO, AM:SAOH, Mo:HCT}}] \label{de:properties}
Let $\ob{L}$ be a complete lattice. Then an operation sequence $(f_n)_{n \in \N}$ on $\ob{L}$ satisfies
\begin{description}
       \item[\HC{1}, boundedness by arguments] if $f_n (\alpha_1,\ldots,\alpha_n) \le \bigmeet_{j=1}^n
       \alpha_j$,
     \item[\HC{2}, monotony] if
       \[ \alpha_1 \le \beta_1,\ldots, \alpha_n \le \beta_n
       \text{ implies }  f_n (\alpha_1,\ldots, \alpha_n) \le f_n
       (\beta_1,\ldots, \beta_n),
       \]
       \item[\HC{3}, omission property] if $f_{n+1} (\alpha_1,\ldots, \alpha_{n+1}) \le f_{n} (\alpha_2, \ldots,
       \alpha_{n+1})$,
       \item[\HC{4}, symmetry] if $f_n (\alpha_1,\ldots, \alpha_n) = f_n (\alpha_{\sigma (1)}, \ldots, \alpha_{\sigma (n)})$ for
                   all $\sigma \in S_n$,
                 \item[\HC{7}, join distributivity] if
                   \begin{multline*} f_n (\alpha_1,\ldots, \alpha_{k-1}, \bigvee_{j \in J} \rho_j, \alpha_{k+1},\ldots, \alpha_n)  =
       \bigvee_{j \in J} f_n (\alpha_1,\ldots, \alpha_{k-1},
       \rho_j, \alpha_{k+1},\ldots, \alpha_n),
                   \end{multline*}
       \item[\HC{8}, nesting property] if $$f_{k}(\alpha_{1},\ldots,\alpha_{k-1},f_{n-k+1}\tup{\alpha}{k}{n})\leq
       f_n\tup{\alpha}{1}{n}$$
  \end{description}
for all $n\in\N$, $k\in\{1,\dots,n\}$, $\alpha_1, \ldots, \alpha_{n+1}, \beta_1,\ldots, \beta_n \in \ob{L}$, and all nonempty families
$(\rho_j)_{j \in J}$ from $\ob{L}$.
 \end{de}
It is known that the commutator operations of many algebras satisfy these properties:
\begin{thm} \label{thm:tcadm}
  Let $\ab{A}$ be an algebra, let $\ob{L}$ its congruence lattice and let $(f_n)_{n \in \N}$ its
  sequence of higher commutator operations as defined in \cite{Bu:OTNO}.
  Then we have
  \begin{enumerate}
  \item The sequence $(f_n)_{n \in \N}$ satisfies \HC1, \HC2, \HC3.
  \item \cite{AM:SAOH} If $\ab{A}$ lies in a congruence permutable variety, then
    $(f_n)_{n \in \N}$ satisfies \HC{4}, \HC{7}, and \HC{8}.
  \item \cite{Mo:HCT} If $\ab{A}$ lies in a congruence modular variety, then
    $(f_n)_{n \in \N}$ satisfies \HC{4} and \HC{7}.
  \end{enumerate}
\end{thm}

This paper is devoted to structures of the form
$\algop{\ob{L}}{\join,\meet, (f_n)_{n \in \N}}$ that
satisfy the properties given in Definition~\ref{de:properties}.
This can be seen as a way to extend the theory of commutator lattices,
which were introduced in \cite{Bi:LT2} and studied in
\cite{Cz:TEDC,Ai:CLFN}.
The structure $\algop{\ob{L}}{\join,\meet, (f_n)_{n \in \N}}$
has been introduced in \cite{AM:SOCO} and was studied further
in \cite{Mu:TLHC}.
In \cite{AM:SOCO}, it was proved that
a finite lattice admits at most $\aleph_0$
many sequences satisfying \HC{3} and \HC{4} (this also follows
from Theorem~\ref{thm:finrep1} below). Very often, such results
on the commutator structure reflect theorems on the clone of
polynomial functions of a universal algebra. For example,
the above result reflects the fact that
on a finite set, there are at most countably many polynomial clones with
a Malcev operation \cite[Theorem~5.3(3)]{Ai:CMCO}. Similar correspondences
can be found between \cite[Theorem~1.2]{AM:SOCO} (there is no infinite
descending chain of commutator operations) and
\cite[Theorem~5.3(1)]{Ai:CMCO} (there is no infinite descending
chain of Mal'cev clones on a finite set),
and between \cite[Theorem~1.1]{AM:SOCO} and \cite[Theorem~1]{AH:CPEO}.
It is worth noting that on the commutator side, it has been proved
that on a finite lattice,
there is no infinite antichain of operation sequences with \HC{3} and
\HC{4} \cite[Theorem~1.2]{AM:SOCO}. The corresponding question on the clone side whether
there is a finite set with an infinite antichain of Mal'cev clones
is still open.

  The goal of this paper is to provide a compact representation
  of the sequence of higher commutator operations
  \[
  C_{\ab{A}} = ([.], \, [.,.], \, [.,.,.], \,\ldots)
  \]
  of a
  finite algebra in a congruence modular variety.
  We will prove that there is a finite set of identities $\Phi$
  of the form
  $[\alpha_1, \ldots, \alpha_n] = \beta$ with $\alpha_1, \ldots,
  \alpha_n, \beta \in \Con (\ab{A})$ such that $C_{\ab{A}}$
  is the largest sequence
  with \HC3 and \HC4 satisfying these identities (Theorem~\ref{thm:largest}).
  It is possible that $C_{\ab{A}}$ is not the only sequence with
  \HC3, \HC4 and $\Phi$. For determining $C_{\ab{A}}$ uniquely,
  Theorem~\ref{thm:unique}
  provides a finite set of \emph{extended commutator equalities},
  which are defined in Section~\ref{sec:commseq}.
  The basic technique is to translate sequences of commutator operations
  into certain order inverting (antitone) functions defined
  on $\N_0^m$.

\section{Translating operation sequences into antitone functions} \label{sec:trans}

When an operation sequence on a finite lattice $\ob{L}$ with $m$ elements
satisfies \HC{4}, then it can
be encoded by a function from $\N_0^m$ to $\ob{L}$.
Suppose that $|\ob{L}| = m$ and
$\ob{L} = \{ \lambda_1, \ldots, \lambda_m\}$.
When $(f_n)_{n \in \N}$ satisfies (HC4),
then using this symmetry,  $f_n (\alpha_1, \ldots, \alpha_n)$
 can be computed as
 \begin{equation} \label{eq:Ff}
 f_n (
     \underbrace{\lambda_1, \ldots, \lambda_1}_{a_1}, \ldots,
     \underbrace{\lambda_m, \ldots, \lambda_m}_{a_m}
     ),
\end{equation}
    where
    $a_j = |\{ k \in \{1, \ldots, m \} \, : \, \alpha_k = \lambda_j\}|$ is the
    number of occurrences of $\lambda_j$ in $(\alpha_1, \ldots, \alpha_n)$.
    For $(a_1, \ldots, a_m) \in \NW{m} \setminus \{(0,\ldots, 0)\}$, let
    $F (a_1, \ldots, a_m)$ be the value of the expression~\eqref{eq:Ff},
    where $n := \sum_{j=1}^m a_j$, and let $F(0, \ldots, 0)$ be the largest
    element of $\ob{L}$.
    We call $F$
    the \emph{encoding} of $(f_n)_{n \in \N}$; $F$ is then a
    function
    from $\NW{m}$ to $\ob{L}$.
    It is easy to see that this encoding provides a bijection
    between the set of operation sequences that satisfy (HC4) and
    the set of all functions from $\NW{m}$ to $\ob{L}$.
    We write elements of $\NW{m}$ in the form
    $\vb{a} = (a_1,\ldots, a_m)$, and we order
    $\NW{m}$ by $\vb{a} \le \vb{b}$ if and only if
    $a_i \le b_i$ for all $i \in \{1, \ldots, m\}$.
    We use $\vb{e}_i$ to the denote the $i$\,th unit vector
    $(0,\ldots, 0,1,0,\ldots, 0)$ with $1$ at the $i$\,th place, and
    we write the supremum of $\vb{a}$ and $\vb{b}$ as
    $\vb{a} \sqcup \vb{b}$ and $\sup \, \{\vb{a}, \vb{b} \}$; it
    can be computed as $(\max(a_1, b_1), \ldots, \max(a_m, b_m))$; dually, we write
    $\vb{a} \sqcap \vb{b}$ for $(\min (a_1, b_1), \ldots, \min(a_m, b_m))$.
    The set $\{1,\ldots, m\}$ will be abbreviated by
    $\ul{m}$, and a function $F : \NW{m} \to \ob{L}$ is called
    \emph{antitone} if $F(\vb{b}) \le F (\vb{a})$ for all
    $\vb{a}, \vb{b} \in \NW{m}$ with $\vb{a} \le \vb{b}$.

    For sequences with \HC{4}, the properties of
    sequences can easily be translated into the following
    properties of their encodings.
\begin{lem} \label{lem:HCencoding}
  Let $\ob{L}$ be a finite lattice with $m := |\ob{L}|$ and
  $\ob{L} = \{ \lambda_1, \ldots, \lambda_m \}$.
  Let $(f_n)_{n \in \N}$ be an operation sequence on $\ob{L}$ that
  satisfies (HC4), and let $F : \NW{m} \to \ob{L}$ be its
  encoding.
  Then $(f_n)_{n \in \N}$ satisfies
  \begin{enumerate}
     \item  \label{it:Ehc1} \HC{1} if and only if for all $\vb{a} \in \NW{m}$ and $j \in \ul{m}$ with $a_j > 0$, we have $F(\vb{a}) \le \lambda_j$.
     \item \label{it:Ehc2} \HC{2} if and only if for all $\vb{a} \in \NW{m}$ and
       $i,j \in \ul{m}$ with $a_j > 0$ and
       $\lambda_i \le \lambda_j$, we have
       $F(\vb{a} - \vb{e}_j + \vb{e}_i) \le F(\vb{a})$.
     \item \label{it:Ehc3} \HC{3} if and only if $F$ is antitone.
     \item \label{it:Ehc7} \HC{7} if and only if for all
       $\vb{a} \in \NW{m}$ and for all $i,j,k \in \ul{m}$ with
       $\lambda_k = \lambda_i \vee \lambda_j$, we have
       $F(\vb{a} + \vb{e}_k) = F(\vb{a} + \vb{e}_i) \join F(\vb{a} + \vb{e}_j)$.
     \item \label{it:Ehc8} \HC{8} if and only if for all $j \in \ul{m}$, for all
       $\vb{a} \in \NW{m}$, for all $\vb{b} \in \NW{m}$
       with $F(\vb{b}) = \lambda_j$ and $\vb{b} \le \vb{a}$,
       we have
       $F (\vb{a} - \vb{b} + \vb{e}_j) \le F (\vb{a})$.
  \end{enumerate}
\end{lem}
\begin{proof}
  For each $n \in \N$ and each $\gb{\alpha} =
  (\alpha_1, \ldots, \alpha_n) \in \ob{L}^n$, we write
  $f (\gb{\alpha})$ for $f_n (\alpha_1, \ldots, \alpha_n)$.
  Furthermore, we define $\vb{a} (\gb{\alpha}) = (a_1 (\gb{\alpha}), \ldots, a_m (\gb{\alpha}))$,
  where for each $i \in \ul{m}$,
  $a_i (\gb{\alpha})$ is defined as the number of occurrences
  of $\lambda_i$ in $\gb{\alpha}$. Formally
  $a_i (\gb{\alpha}) = \left|\{ j \in \ul{n} \,\,\colon\, \alpha_j = \lambda_i \}\right|$.
  For $\vb{a} \in \NW{m}$, we define $\gb{\lambda} (\vb{a}) \in \ob{L}^{\sum_{j=1}^m a_j}$
  by
  \[
       \gb{\lambda} (\vb{a}) = (\underbrace{\lambda_1, \ldots, \lambda_1}_{a_1 \text{ times}}, \ldots,
       \underbrace{\lambda_m, \ldots, \lambda_m}_{a_m \text{ times}}).
       \]
       With this notation, the definition of $F(\vb{a})$, for $\vb{a} \in \NW{m} \setminus
       \{ (0,\ldots,0)\}$, can be stated
  as
  \begin{equation} \label{eq:encodingeq}
    F (\vb{a}) = f  (\gb{\lambda} (\vb{a})).
  \end{equation}
  Then the symmetry property \HC{4} yields that
  \begin{multline} \label{eq:symm}
    f (\gb{\alpha}) = f_n (\alpha_1, \ldots, \alpha_n)
    \\ =
    f_n (\underbrace{\lambda_1, \ldots, \lambda_1}_{a_1 (\gb{\alpha}) \text{ times}},
         \ldots,
         \underbrace{\lambda_m, \ldots, \lambda_m}_{a_m (\gb{\alpha}) \text{ times}})
    =
     f (\gb{\lambda} (\vb{a} (\gb{\alpha})))
        = F (\vb{a} (\gb{\alpha})).
   \end{multline}
  To prove \eqref{it:Ehc1}, we first assume that $(f_n)_{n \in \N}$ satisfies \HC{1}.
  Let $j \in \ul{m}$ and let $\vb{a} \in \NW{m}$ with $a_j > 0$.
  Then $\lambda_j$ appears $a_j$ times in $\gb{\lambda} (\vb{a})$, and
 therefore~\eqref{eq:encodingeq} and \HC{1} yield $F(\vb{a}) \le \lambda_j$.
  For the ``if''-direction of item~\eqref{it:Ehc1}, we
  fix $n \in \N$ and $\gb{\alpha} \in \ob{L}^n$.
  Then by~\eqref{eq:symm},
  $f (\gb{\alpha} ) = F (\vb{a} (\gb{\alpha}))$. Let $i \in \ul{n}$, and
  let $j$ be such that $\lambda_j = \alpha_i$. Then
  $a_j (\gb{\alpha}) > 0$, and therefore by the assumption,
  $F (\vb{a} (\gb{\alpha})) \le \lambda_j = \alpha_i$.
  Combining this for all $i \in \ul{n}$, we obtain
  $F (\vb{a} (\gb{\alpha})) \le \bigmeet_{i \in \ul{n}} \alpha_i$,
  implying $f (\gb{\alpha} ) \le \bigmeet_{i \in \ul{n}} \alpha_i$.

  For proving \eqref{it:Ehc2}, we assume that $(f_n)_{n \in \N}$ satisfies
  \HC{2}. Let $\vb{a} \in \NW{m}$ and
  $i,j \in \ul{m}$ be such that $a_j > 0$ and
  $\lambda_i \le \lambda_j$.
  Then $F(\vb{a}) = f (\gb{\lambda} (\vb{a}))$, which is
  equal to
  $f (\lambda_j, \gb{\lambda} (\vb{a} - \vb{e}_j))$ by \HC4. Using \HC{2},
  we have
  $f (\lambda_j, \gb{\lambda} (\vb{a} - \vb{e}_j)) \ge
  f (\lambda_i, \gb{\lambda} (\vb{a} - \vb{e}_j))$. By \HC{4}, the last
  expression is equal to
  $f (\gb{\lambda} (\vb{a} - \vb{e}_j + \vb{e}_i)) =
   F (\vb{a} - \vb{e}_j + \vb{e}_j)$.
  For the ``if''-direction of~\eqref{it:Ehc2}, we notice that by \HC4, it is sufficient
  to show that for all $\gb{\alpha} \in \ob{L}^n$ and $\beta_1 \in \ob{L}$ with
  $\alpha_1 \le \beta_1$, we have
  \begin{equation} \label{eq:mon1}
    f_n (\alpha_1, \alpha_2, \ldots, \alpha_n) \le
    f_n (\beta_1, \alpha_2, \ldots, \alpha_n).
  \end{equation}
  To prove~\eqref{eq:mon1}, we choose $i,j \in \ul{m}$ such that
  $\alpha_1 = \lambda_i$ and $\beta_1 = \lambda_j$.
  It is then easy to see that
  $\vb{a} := \vb{a} (\beta_1, \alpha_2, \ldots, \alpha_n)$ satisfies
  $a_j > 0$.
  Since $\gb{\alpha}$ is obtained from
  $(\beta_1, \alpha_2, \ldots, \alpha_n)$ by
  replacing the leading $\beta_1 =\lambda_j$ with
  $\alpha_1 = \lambda_i$, we have $\vb{a} (\gb{\alpha}) = \vb{a} - \vb{e}_j + \vb{e}_i$.
  Now by \eqref{eq:symm}, we have
  $f_n (\alpha_1, \alpha_2, \ldots, \alpha_n)
  =
  F (\vb{a} (\gb{\alpha}))  =
  F (\vb{a} - \vb{e}_j + \vb{e}_i)$. By assumption,
  $F (\vb{a} - \vb{e}_j + \vb{e}_i) \le F (\vb{a}) =
  F (\vb{a} (\beta_1, \alpha_2, \ldots, \alpha_n))$, which is equal to
  $f_n (\beta_1, \alpha_2, \ldots, \alpha_n)$ by
  \eqref{eq:symm}. This establishes~\eqref{eq:mon1}.

   For proving~\eqref{it:Ehc3}, we notice that
  $F$ is antitone if and only if for all $\vb{a} \in \NW{m}$ and
  $i \in \ul{m}$, we have
  \begin{equation} \label{eq:simpleanti}
    F(\vb{a} + \vb{e}_i) \le F (\vb{a}).
  \end{equation}
  We
  first assume that $(f_n)_{n \in \N}$ satisfies
  \HC3.
  In the case $\vb{a} = (0,\ldots,0)$, \eqref{eq:simpleanti} holds
  by the definition of $F(0,\ldots,0)$ as the maximal element of $\ob{L}$.
  If $\vb{a} \neq (0,\ldots,0)$, we have (using~\eqref{eq:encodingeq}) that
  $F(\vb{a} + \vb{e}_i) =
  f (\gb{\lambda} (\vb{a} + \vb{e}_i))$, which by \HC4 is equal to
  $f (\lambda_i, \gb{\lambda} (\vb{a}))$. Now by \HC3, the last expression
  is $\le f (\gb{\lambda} (\vb{a})) = F (\vb{a})$. This completes the proof
  of \eqref{eq:simpleanti}.
  For the ``if''-direction of \eqref{it:Ehc3}, we let $n \in \N$ and
  $(\alpha_1, \ldots, \alpha_{n+1}) \in \ob{L}^{n+1}$.
  Let $i \in \ul{m}$ be such that $\lambda_i = \alpha_1$. Then
  we have
  $f_{n+1} (\alpha_1, \ldots, \alpha_{n+1}) =
    F (\vb{a} (\alpha_1, \ldots, \alpha_{n+1}))
    = F (\vb{a} (\alpha_2, \ldots, \alpha_{n+1}) + \vb{e}_i)$.
    By assumption, the last expression is $\le
    F (\vb{a} (\alpha_2, \ldots, \alpha_{n+1})) =
    f_n (\alpha_2, \ldots, \alpha_{n+1}).$

    For proving~\eqref{it:Ehc7}, we notice that since $\ob{L}$ is finite
    and $(f_n)_{n \in \N}$ satisfies \HC4, the property \HC7 is equivalent
    to
    \[
    f_n (\lambda_i \join \lambda_j, \alpha_2, \ldots, \alpha_n)
    =
    f_n (\lambda_i, \alpha_2, \ldots, \alpha_n)
     \join
     f_n (\lambda_j, \alpha_2, \ldots, \alpha_n)
     \]
     for all $n \in \N$, $i,j \in \ul{m}$ and
     $\alpha_2, \ldots, \alpha_n$.
     Let $k \in \ul{m}$ be such that
     $\lambda_k = \lambda_i \join \lambda_j$.
     Let $\vb{a} := \vb{a} (\alpha_2, \ldots, \alpha_n)$.
     Then $\vb{a} (\lambda_i \join \lambda_j, \alpha_2, \ldots, \alpha_n)
           = \vb{a} (\lambda_k, \alpha_2, \ldots, \alpha_n)
           = \vb{a} + \vb{e}_k$,
     $\vb{a} (\lambda_i , \alpha_2, \ldots, \alpha_n) = \vb{a} + \vb{e}_i$, and
           $\vb{a} (\lambda_j , \alpha_2, \ldots, \alpha_n) = \vb{a} + \vb{e}_j$.
           From this, one easily infers that \HC7 is equivalent to
           the condition stated in item~\eqref{it:Ehc7}.

           For~\eqref{it:Ehc8}, we observe that if $\lambda_j = f_{n-k+1} (\alpha_k, \ldots, \alpha_n)$
           then
           for $\vb{a} := \vb{a} (\alpha_1, \ldots, \alpha_n)$ and
           $\vb{b} := \vb{a} (\alpha_k, \ldots, \alpha_n)$, we have
           $\vb{b} \le \vb{a}$,
           $\vb{a} (\alpha_1, \ldots, \alpha_{k-1}, \lambda_j) =
           \vb{a} - \vb{b} + \vb{e}_j$,
           which implies
           $f_k (\alpha_1, \ldots, \alpha_{k-1}, \lambda_j) =
            F (\vb{a} - \vb{b} + \vb{e}_j)$, and
            $f_n (\alpha_1, \ldots, \alpha_n) = F (\vb{a})$.
              From this, one easily infers that \HC8 is equivalent to
           the condition stated in item~\eqref{it:Ehc8}.
\end{proof}
For a lattice $\ob{L} = \{\lambda_1, \ldots, \lambda_m\}$
with $m$ elements and $F : \NW{m} \to \ob{L}$,
we say that $F$ satisfies a condition \HC{x} (with
$x \in \{1,2,3,7,8\}$) when the corresponding condition given in
Lemma~\ref{lem:HCencoding} is satisfied.
The set of all operation sequences on $\ob{L}$ is isomorphic
to the (complete) lattice $\ob{L}^{\bigcup_{i \in \N} \ob{L}^i}$ via the isomorphism
$\Phi$ given by
$\Phi ((f_i)_{i \in \N}) \,\, (\alpha_1, \ldots, \alpha_n)
=
f_n (\alpha_1, \ldots, \alpha_n)$.
The operation sequences satisfying (HC4) then correspond to
a sublattice which, by the encoding given above, is
isomorphic to $\ob{L}^{\N_0^m}$.
In this note, we will mainly study those operation sequences
that satisfy (HC3) and (HC4). These sequences form a complete sublattice
of $\ob{L}^{\bigcup_{i \in \N} \ob{L}^i}$ \cite[Proposition~2.6]{Mu:TLHC}, which, by
Lemma~\ref{lem:HCencoding}\eqref{it:Ehc3},
is isomorphic
to the lattice of antitone function from $(\NW{m}, \le)$ to
$\ob{L}$. More information on these lattices is contained in
\cite{Mu:TLHC}.

\section{Representing antitone functions}
In this section, we will provide some finite representations
of an antitone function $F$ from $\NW{m}$ into a finite lattice
$\ob{L}$. We will identify the function $F$ with its graph
$\{ (\vb{a}, F(\vb{a})) \mid \vb{a} \in \NW{m} \}$, which is a
subset of $\NW{m} \times \ob{L}$.
First, we seek to interpolate a subset $G$ of $\NW{m} \times \ob{L}$
by an antitone function. For a complete lattice $\ob{L}$ and
a subset $S$ of $\ob{L}$, the meet (or infimum) of all elements in $S$
will be denoted by $\bigmeet S$, $\bigmeet_{\sigma \in S} \sigma$ or
$\bigmeet_{\lambda \in \ob{L}, \lambda \in S} \lambda$, and
we define  the empty intersection $\bigmeet \emptyset$ as the largest element of $\ob{L}$.

\begin{de} \label{de:fga}
  Let $\ob{L}$ a complete lattice, and
  let $G \subseteq \NW{m} \times \ob{L}$. Then the
  function $F_G$ from $\NW{m}$ to $\ob{L}$ \emph{represented by $G$}
  is defined by
  \begin{equation} \label{eq:fga}
  F_G (\vb{x}) := \bigmeet \{ \gamma \mid
  (\vb{c}, \gamma) \in G \text{ and }
  \vb{c} \le \vb{x} \}
  \end{equation}
  for all $\vb{x} \in \NW{m}$.
\end{de}
\begin{lem} \label{lem:largest}
    Let $\ob{L}$ be a complete lattice, and
    let $G \subseteq \NW{m} \times \ob{L}$.
    Then
    \begin{enumerate}
       \item \label{it:l1} $F_G$ is antitone, i.e., for all $\vb{a}, \vb{b} \in \NW{m}$ with
       $\vb{a} \le \vb{b}$, we have $F_G(\vb{b}) \le F_G(\vb{a})$.
       \item \label{it:l2} $F_G$ is the largest antitone function $F$ from
         $\NW{m} \to \ob{L}$ such that
         \begin{equation} \label{eq:FG}
             \text{for all }
             (\vb{c}, \gamma) \in G, \text{ we have } F (\vb{c}) \le \gamma.
         \end{equation}
       \item  \label{it:l3}
          Suppose that
    for all $(\vb{a}, \alpha)$ and $(\vb{b}, \beta) \in G$
    with $\vb{a} \le \vb{b}$, we have
    $\beta \le \alpha$.
    Then $G \subseteq F_G$.
     \end{enumerate}
  \end{lem}
  \begin{proof}
  \eqref{it:l1}
  Since $\vb{a} \le \vb{b}$, we have
  \( \{ (\vb{c}, \gamma) \in G \mid \gamma \le \vb{a} \} \subseteq
  \{ (\vb{c}, \gamma) \in G \mid \gamma \le \vb{b} \} \)
  and therefore $F_G (\vb{b}) \le F_G (\vb{a})$.

  \eqref{it:l2}
  We first show that $F_G$ satisfies~\eqref{eq:FG}.
  If $(\vb{c}, \gamma) \in G$, then $\gamma$ is one of the
  elements of $\ob{L}$ that appear in the intersection defining
  $F_G (\vb{c})$. Thus $F_G (\vb{c}) \le \gamma$, and
  therefore $F_G$ satisfies~\eqref{eq:FG}.
  To show that $F_G$ is the largest antitone function
  with \eqref{eq:FG},
  we let
  $F$ be an antitone function satisfying~\eqref{eq:FG}, and
  we let $\vb{a} \in \NW{m}$. Then for each $(\vb{c}, \gamma) \in G$
  with $\vb{c} \le \vb{a}$, we have
  $F(\vb{a}) \le F (\vb{c})$ because
  $F$ is antitone, and $F(\vb{c}) \le \gamma$ because
  of the assumptions on $F$, and thus
  $F(\vb{a}) \le \gamma$. Therefore,
  $F(\vb{a}) \le  \bigmeet \{ \gamma \mid
  (\vb{c}, \gamma) \in G \text{ and }
  \vb{c} \le \vb{a} \} = F_G (\vb{a})$.

  \eqref{it:l3}
    Let $(\vb{b}, \beta) \in G$. We compute
    $F_G (\vb{b})$ as $\bigmeet A$ with $A =
\{ \alpha \mid
  (\vb{a}, \alpha) \in G \text{ and }
\vb{a} \le \vb{b} \}$.
Since $(\vb{b}, \beta) \in G$ and $\vb{b} \le \vb{b}$,
we have $\beta \in A$. Furthermore, by the assumption on $G$,
all $\alpha \in A$ satisfy
$\beta \le \alpha$. This implies $F_G (\vb{b}) = \beta$.
\end{proof}
  For a subset $A$ of a partially ordered set $(W, \le)$, the set
  of \emph{minimal elements} of $A$,
  $\Min (A)$, is given by
  \[
   \Min (A) =
  \{ a \in A \mid \text{ for all } b \in A \,:\,
     b \le a \Rightarrow b = a \}.
  \]
  The set of \emph{maximal elements} of $A$, $\Max (A)$, is defined
  dually.
  $(W, \le)$ is called \emph{well partially ordered}
  if $(W, \le)$ has no infinite descending chains and
  and no infinite antichains.
  For $m \in \N$, the product
  $W^m$ is ordered by $(w_1, \ldots, w_m) \le (w'_1, \ldots, w'_m)$ if
  $w_i \le w'_i$ for all $i \in \ul{m}$. We will need
  the fact that being well partially ordered is preserved under
  forming finite products; this is well-known and can be found,
  e.g., in \cite{AH:FTIS}.
  \begin{lem} \label{lem:dick} \cite[p.195, Example~(4)]{AH:FTIS}
      Let $m \in \N$.
   If  $(W, \le)$ is a well partially ordered set, then also
   $(W^m, \le)$ is well partially ordered.
   In particular, $(\NW{m}, \le)$ has no infinite antichains.
  \end{lem}
  The fact that $\NW{m}$ has no infinite antichains is often
  called Dickson's Lemma \cite{Di:FOTO}, and a more detailed
  discussion with proofs is contained, e.g.,  in the expository
  survey \cite{AA:DLHT}.
  As a consequence, for every subset $A$ of $\NW{m}$, the
  set $\Min (A)$ of minimal elements of $A$ is finite.
  Furthermore, since there are no infinite descending chains in $A$,
  we obtain that for every
  $\vb{a} \in A$, there is a $\vb{b} \in \Min (A)$ with
  $\vb{b} \le \vb{a}$.
 \begin{thm} \label{thm:finrep1}
  Let $\ob{L}$ be a finite lattice. Then every antitone function $F$
  from $\NW{m}$ to $\ob{L}$ can be represented by a finite subset $G$
  of $F$.
\end{thm}
\begin{proof}
For every $\alpha \in \ob{L}$, the set $F^{-1} (\{\alpha\})$ contains
only finitely many minimal elements. Let
\begin{multline*}
  \Min (F^{-1} (\{\alpha\})) \\ = \{ \vb{a} \in \NW{m} \mid F(\vb{a}) = \alpha \text{ and }
\neg \exists \vb{b} \in  \NW{m} \,  :   \,
(\vb{b} < \vb{a} \text{ and } F (\vb{b}) = \alpha)  \}
\end{multline*}
be the set of these minimal elements.
We note that this set is empty if $\alpha$ is not an element
of the image of $F$. Let
$G := \bigcup_{\alpha \in \ob{L}}
(\Min (F^{-1} (\{\alpha\})) \times \{ \alpha \})$; this set $G$ is
a subset of $F$.
 We claim that $F_G = F$.
 To show this, we fix $\vb{a} \in \NW{m}$ and compute
  $F_G(\vb{a}) =  \bigmeet B$, where $B := \{ \beta \mid
  (\vb{b}, \beta) \in G \text{ and }
  \vb{b} \le \vb{a} \}$.
   Let $\alpha := F (\vb{a})$. Then there exists
   $\vb{b} \in \Min (F^{-1} (\{\alpha\}))$ with
   $\vb{b} \le \vb{a}$. Then  $(\vb{b}, \alpha) \in G$,
   and therefore $\alpha \in B$.
   Thus $F_G (\vb{a}) \le F (\vb{a})$.
  For proving $F(\vb{a}) \le F_G (\vb{a})$, we observe that
  every $\beta \in B$ satisfies $\alpha \le \beta$; to see this,
  we fix $(\vb{b}, \beta) \in G$ with $\vb{b} \le \vb{a}$.
  Since $F(\vb{b}) = \beta$, the fact that
  $F$ is antitone implies $\alpha \le \beta$. Hence
  $\alpha \le \bigmeet B = F_G (\vb{a})$.
\end{proof}
The representation constructed in this proof will be
called the \emph{canonical representation} of $F$:
\begin{de}
Let $m \in \N$, let $\ob{L}$ be a finite lattice,
and let $F$ be an antitone function from
$\NW{m}$ to $\ob{L}$.
Then the \emph{canonical representation} $G$ of $F$ is defined
by
$G := \bigcup_{\alpha \in \ob{L}}
\{ (\vb{a}, \alpha) \mid \alpha \in \ob{L} \text{ and }
\vb{a} \text{ is minimal in }\NW{m}\text{ with }
     F (\vb{a}) = \alpha \}$.
\end{de}

For a lattice $\ob{L}$ with maximal element $1$,
$\beta \in \ob{L}$ and $\vb{b} \in \NW{m}$, we let
$\one{\vb{b}}{\beta}$ be the function from $\NW{m}$ to $\ob{L}$ defined
by $\one{\vb{b}}{\beta} (\vb{x}) = \beta$ if $\vb{x} \ge \vb{b}$ and
$\one{\vb{b}}{\beta} (\vb{x}) = 1$ if $\vb{x} \not\ge \vb{b}$.
\begin{cor}
  Let $\ob{L}$ be a finite lattice, and let
  $f : \NW{m} \to \ob{L}$ be an antitone function. Then there
  is a finite subset $G$ of $\NW{m} \times \ob{L}$ such that
  \begin{equation} \label{eq:fisfg}
  f (\vb{x}) = \bigmeet_{(\vb{b}, \beta) \in G} \one{\vb{b}}{\beta} (\vb{x})
  \end{equation}
  for all $\vb{x} \in \NW{m}$. In other words, the semilattice of
  antitone functions from $\NW{m}$ to the semilattice $(\ob{L}, \meet)$
  is generated by $\{\one{\vb{b}}{\beta} \mid
  \vb{b} \in \NW{m}, \beta \in \ob{L} \}$.
\end{cor}
\begin{proof}
  Theorem~\ref{thm:finrep1} provides us with a finite set $G$
  such that
  \begin{multline*}
     f (\vb{x}) = \bigmeet \{ \beta \mid (\vb{b}, \beta) \in G, \vb{b} \le \vb{x} \}
  \\ = (\bigmeet \{ \beta \mid (\vb{b}, \beta) \in G, \vb{b} \le \vb{x} \})
  \meet
    (\bigmeet\{ 1 \mid (\vb{b}, \beta) \in G, \vb{b} \not\le \vb{x} \})
  \\ =
  \bigmeet \{ \one{\vb{b}}{\beta} (\vb{x}) \mid (\vb{b}, \beta) \in G \},
  \end{multline*}
  and therefore~\eqref{eq:fisfg} is satisfied.
  From this equation, we see that $f$ is the meet
  of $|G|$ many functions $\one{\vb{b}}{\beta}$
  with  $\vb{b} \in \NW{m}$ and  $\beta \in \ob{L}$.
  If $G \neq \emptyset$, this shows that $f$ lies in the
  semilattice generated by the functions $\one{\vb{b}}{\beta}$.
  In the case $G = \emptyset$, we have $f(\vb{x}) = 1$
  for all $\vb{x} \in \NW{m}$, and thus
  $f = \one{(0,\ldots, 0)}{1}$, where $1$ denotes the largest element
  of $\ob{L}$.
\end{proof}

It may happen that
the antitone function $F$ is represented by $G$,
but there exists $F_1 \neq F$ with $G \subseteq F_1$.
 The following proposition characterizes
which antitone functions can be uniquely determined by a finite
subset of their graph.
\begin{pro} \label{pro:nonunique}
  Let $\ob{L}$ be a finite lattice with minimal element $0$,
  and let $F : \NW{m} \to \ob{L}$ be an antitone function.
  Then the following are equivalent:
  \begin{enumerate}
  \item \label{it:f1}  There is a finite subset $G$ of $F$ such that $F$ is
    the unique antitone function with $G \subseteq F$.
  \item \label{it:f2} For each $i \in \ul{m}$ there is $N \in \N$ such
    that for all $c \in \N$ with $c \ge N$, we have
    $F (c \vb{e}_i) = 0$.
  \end{enumerate}
\end{pro}
\begin{proof}
  \eqref{it:f1}$\Rightarrow$\eqref{it:f2}:
  We assume that there is an $i \in \ul{m}$ such that
  $F(c \vb{e_i}) \neq 0$ for infinitely many $c \in \N$.
Since $F$ is antitone and $\ob{L}$ is finite, the sequence $(F (c
\vb{e}_i))_{c \in \N}$ must eventually become constant; let $N \in
\N$ and $\gamma \in \ob{L}$ be such that $F (c \vb{e}_i) = \gamma$
for all $c \ge N$. Since $F(c \vb{e_i}) \neq 0$ for infinitely many
$c \in \N$, we have $\gamma > 0$. Seeking to prove the negation
of~\eqref{it:f1}, we prove that for every finite subset $G$ of $F$,
there are at least two antitone functions containing $G$ as a
subset. Since $G$ is finite, there is $c \in \N$ with $c \ge N$, and
for all $(\vb{a}, \alpha) \in G$, $c  \vb{e}_i \not\le \vb{a}$. Now
let $G' := G \cup \{ (c  \vb{e_i}, 0) \}$. We will now show that
$G'$ satisfies the assumptions of
Lemma~\ref{lem:largest}\eqref{it:l3}. To this end, let $(\vb{a},
\alpha)$ and $(\vb{b}, \beta) \in G'$ with $\vb{b} \le \vb{a}$. If
$(\vb{a}, \alpha)$ and $(\vb{b}, \beta)$ are both elements of $G$,
then $F(\vb{a}) \le F(\vb{b})$ and therefore $\alpha \le \beta$. If
$\vb{a} = c \vb{e}_i$, then $0 \le \beta$ because $0$ is the minimal
element of $\ob{L}$. Finally, if $\vb{b} = c \vb{e}_i$, then the
elements $(\vb{a}, \alpha)$ from $G$ do not satisfy $\vb{b} \le
\vb{a}$, and therefore $(\vb{a}, \alpha) \in G' \setminus G$. Hence
$(\vb{a}, \alpha) = (c \vb{e}_i, 0)$ and thus $\alpha = \beta = 0$,
implying that $\alpha \le \beta$ also holds in this case. Therefore
the assumptions of Lemma~\ref{lem:largest}\eqref{it:l3} are
satisfied, and we have $G'  \subseteq F_{G'}$. Now $F$ and $F_{G'}$
are two antitone functions containing $G$, but $F (c \vb{e}_i) =
\gamma$ and $F_{G'} (c \vb{e}_i) = 0$.

\eqref{it:f2}$\Rightarrow$\eqref{it:f1}:
We first show that for every $\alpha \in \ob{L}$ with $\alpha > 0$, the
set $\{ \vb{a} \mid F (\vb{a}) = \alpha \}$ is finite.
Suppose $\{ \vb{a} \mid F (\vb{a}) = \alpha \}$ is infinite.
Then there is $i \in \ul{m}$ such that
$B := \{ a_i \mid F (a_1, \ldots, a_m) = \alpha \}$ is infinite.
Then for every $c \in \N_0$, we have an $\vb{a} \in F^{-1} (\{\alpha\})$
with  $c \vb{e}_i \le \vb{a}$, and therefore
$F (c \vb{e}_i) \ge F (\vb{a}) = \alpha$, contradicting the assumption.
Hence the set $\{ \vb{a} \mid F (\vb{a}) = \alpha \} = F^{-1} (\{ \alpha \})$
is finite, and therefore the
set of its maximal elements $\Max (F^{-1} ( \alpha ))$ is finite.
Now let
\[
G := \bigcup_{\alpha \in \ob{L}} \Min (F^{-1} (\alpha)) \times \{\alpha\}
\,\,\, \cup \,\,\,
\bigcup_{\alpha \in \ob{L}\setminus \{0\} } \Max (F^{-1} (\alpha)) \times \{\alpha\}.
\]
Let $H$ be an antitone function with $G \subseteq H$. We show $H = F$. To this
end, we fix $\vb{a} \in \NW{m}$ and set $\alpha := F (\vb{a})$.
Then there is a minimal element $\vb{b} \in F^{-1} (\{\alpha\})$
with $\vb{b} \le \vb{a}$. If $\alpha = 0$, then
$H (\vb{a}) \le H (\vb{b}) = F(\vb{b}) = 0$.
If $\alpha > 0$, there is
a maximal element $\vb{c} \in F^{-1} (\{\alpha\})$ with
$\vb{b} \le \vb{a} \le \vb{c}$.
Then $\alpha = H (\vb{b}) \ge H (\vb{a}) \ge H(\vb{c}) = \alpha$, which implies
$H(\vb{a}) = \alpha$.
\end{proof}

\section{Extensions of antitone functions}
From Proposition~\ref{pro:nonunique}, we gather
that many antitone functions from $\NW{m}$ to a
finite lattice $\ob{L}$ are never completely determined by a finite subset
of their graph. However, if we extend the domain from $\NW{m}$ to $\NWC{m}$,
then such a finite representation is possible.
Here, denoting the ordinal number $\N_0$ by $\omega$, we see
$\N_0 \cup \{\infty\}$ as the ordinal number $\omega^+$, the successor
of $\omega$. Therefore the ordering of $\N_0 \cup \{\infty\}$ has
$\infty$ as its largest element, and we understand $\NWC{m}$ to be ordered
by the product order.
For an antitone $F$ from $\NW{m}$ into the complete lattice $\ob{L}$, we define its extension
$\ext{F} : \NWC{m} \to \ob{L}$
by
\begin{equation} \label{eq:Fhat}
   \ext{F} (\vb{x}) := \bigmeet \{ F (\vb{b}) \mid \vb{b} \le \vb{x},
   \vb{b} \in \NWO{m} \}
\end{equation}
   for all $\vb{x} \in \NWC{m}$.
From its definition, we see that $\ext{F}$ is antitone and agrees with
$F$ on $\NW{m}$. However, not every antitone $H : \NWC{m} \to \ob{L}$
is the extension of its restriction to $\NW{m}$. As an example,
we set $H (\vb{a} ) := 1$ if $\vb{a} \in \NW{m}$ and
$H(\vb{a}) := 0$ if one of the components of $\vb{a}$ is $\infty$, where
$0$ and $1$ are the smallest and largest element of $\ob{L}$, respectively.
Then $H$ is antitone, but the extension of $H|_{\NW{m}}$ is the constant
function with value~$1$. With suitable topologies on $\NWC{m}$ and $\ob{L}$
(on $\NWC{m}$, we take the product topology of the Alexandroff extension of the discrete topology on $\N_0$,
and on $\ob{L}$ we take the discrete topology), we see that this function $H$ is not
continuous, whereas every extension $\ext{F}$ is continuous. We will not pursue or make use of
this topological view, but we consider it worth noting that we
succeed in representing the \emph{continuous} antitone functions by finite subsets of their
graphs.
The following lemma tells that the extension of $F_G$ can be evaluated
in the same way as $F_G$.
\begin{lem} \label{lem:fga}
  Let $\ob{L}$ be a finite lattice, and let
  $G$ be a finite subset of $\NW{m} \times \ob{L}$.
  Then we have:
  \begin{enumerate}
     \item \label{it:fg1} For every $\vb{a} \in \NWC{m}$,
  \(
      \ext{F_G} (\vb{a}) = \bigmeet \{ \gamma \mid
  (\vb{c}, \gamma) \in G \text{ and }
      \vb{c} \le \vb{a} \}.
      \)
    \item \label{it:fg15}
      For every $\vb{a} \in \NWC{m}$, there is
      $\vb{b} \in \NW{m}$ with $\vb{b} \le \vb{a}$ and
      $F_G (\vb{b}) = \ext{F_G} (\vb{a})$.
\item \label{it:fg2}
  Suppose that $\alpha \in \ob{L}$ and $\vb{a} \in \NWC{m}$
  is minimal with $\ext{F_G} (\vb{a}) \le \alpha$. Then
  $\vb{a} \in \NW{m}$.
 \end{enumerate}
\end{lem}
  \begin{proof}
    \eqref{it:fg1}
  We have
  \begin{multline*}
     \ext{F_G} (\vb{a})  =  \bigmeet \{ F_G (\vb{b}) \mid \vb{b} \le \vb{a},
     \vb{b} \in \NW{m} \}
     \\ =
     \bigmeet \{
     \bigmeet_{(\vb{c}, \gamma) \in G, \vb{c} \le \vb{b}} \gamma \mid
     \vb{b} \le \vb{a},
     \vb{b} \in \NW{m} \} =
     \bigmeet_{(\vb{c}, \gamma) \in G, \vb{c} \le \vb{a}} \gamma.
  \end{multline*}

   \eqref{it:fg15}
  For each $\vb{c} \in \NWC{m}$, we define the set $B (\vb{c})$
  by
  \[
  B (\vb{c}) := \{ (\vb{d}, \delta) \in G \mid \vb{d} \le \vb{c} \}.
  \]
  Then $\vb{b} := \sup \, \{ \vb{d} \mid (\vb{d}, \delta) \in B (\vb{a}) \}$
  satisfies $\vb{b} \le \vb{a}$.
  Furthermore, by the construction of $\vb{b}$,
  we have $B (\vb{a}) \subseteq B(\vb{b})$.
     Using item~\eqref{it:fg1} to compute the values of  $\ext{F_G}$, we obtain
     $\ext{F_G} (\vb{b}) \le \ext{F_G} (\vb{a})$. The converse
     $\ext{F_G} (\vb{a}) \le \ext{F_G} (\vb{b})$ follows from
     the antitony of $\ext{F_G}$, and therefore
     $\ext{F_G} (\vb{b}) = \ext{F_G} (\vb{a})$.
      Since $G$ is a finite subset of $\NW{m} \times \ob{L}$,
      $\vb{b}$ is the supremum of finitely many elements of $\NW{m}$,
      and therefore $\vb{b} \in \NW{m}$.
      Finally, we observe that $F_G (\vb{b}) = \ext{F_G} (\vb{b})$;
      this is a consequence of item~\eqref{it:fg1}, or can be
      seen directly from the definition of the extension $\ext{F_G}$;
      altogether, we have $\vb{b} \in \NW{m}$ and
      $F_G (\vb{b}) = \ext{F_G} (\vb{a})$.

\eqref{it:fg2}
    follows directly from~\eqref{it:fg15}.
\end{proof}

  For an element $\vb{a}$ of a partially ordered set $\ob{P}$,
  we write $\vb{a}{\uparrow}$ for the upward closed set
  $\{ \vb{b} \in \ob{P} \mid \vb{a} \le \vb{b} \}$, and
  for $A \subseteq \ob{P}$, we let
  $A{\uparrow} := \bigcup \{ \vb{a}{\uparrow} \,\mid\, \vb{a} \in A \}$
  be the smallest upward closed set that contains $A$ as a subset.
\begin{pro} \label{pro:max}
   Let
   $A$ be a finite subset of $\NW{m}$, and let
   $\vb{b} \in \NWC{m} \setminus (A{\uparrow})$.
   Then
   $C := \{\vb{c} \in \NWC{m} \setminus (A{\uparrow})  \mid
            \vb{b} \le  \vb{c} \}
            =  \{ \vb{c} \in \NWC{m} \mid \vb{b} \le  \vb{c}, \forall \vb{a} \in
          A \, : \, \vb{a} \not\le \vb{c} \}$
   contains a maximal element.
\end{pro}
\begin{proof}
  We use Zorn's Lemma to find a maximal element
  of $C$. To this end, let $S$ be a nonempty linearly ordered
  subset of $C$, and let $\vb{s} \in \NWC{m}$ be the
  supremum of $S$. We want to show that $\vb{s}$ is an element of
  $\NWC{m} \setminus (A{\uparrow})$. Seeking a contradiction,
  we suppose that there is $\vb{a} \in A$ such that
  $\vb{s} \ge \vb{a}$.
  For every $i \in \ul{m}$,
  $a_i \le \sup \{ c_i \mid \vb{c} \in S\}$. Since
  $a_i \in \N_0$, there is $\vb{c}^{(i)} \in S$ with
  $a_i \le (\vb{c}^{(i)})_i$. Hence
  for
  \[
  \vb{d} := \vb{c}^{(1)} \sqcup \cdots \sqcup \vb{c}^{(m)},
  \]
  we have
  $\vb{a} \le \vb{d}$ and therefore
  $\vb{d} \in A{\uparrow}$. Since $S$ is linearly ordered,
  there is $i \in \ul{m}$ such that $\vb{d} = \vb{c}^{(i)}$,
  and therefore
  $\vb{d} \in S$.
  Since $S \subseteq C$, we then have $\vb{d} \in C$ and thus
  $\vb{d} \not\in A{\uparrow}$.
  This contradiction establishes
  $\vb{s} \in \NWC{m} \setminus (A{\uparrow})$.
  Since $\vb{s} \ge \vb{b}$ by its construction, we have
  $\vb{s} \in C$, and thus
  the set $S$ has $\vb{s}$ as an upper bound in $C$.
  Now Zorn's Lemma implies the existence of a maximal element
  in $C$.
\end{proof}

We will now see that many antitone functions are completely determined
by their values on a finite subset of $\NWC{m}$.
\begin{lem} \label{lem:unique}
  Let $m \in \N$, let $\ob{L}$ be a complete lattice, and let
  $F : \NW{m} \to \ob{L}$.
   For every $\alpha \in \ob{L}$, let
  $P_{\alpha}$ be the set of all $\vb{a} \in \NWC{m}$
  that are minimal with $\ext{F} (\vb{a}) \le \alpha$ or
  maximal with $\ext{F} (\vb{a}) \not\le \alpha$, and let
  $P := \bigcup_{\alpha \in \ob{L}} P_{\alpha}$.
  Then for every antitone function $F_1 : \NWC{m} \to \ob{L}$
  with $F_1|_{P} = \ext{F}|_P$, we have $F_1 = \ext{F}$.
\end{lem}
\begin{proof}
   Let $\vb{x} \in \NWC{m}$, and let $\beta := \ext{F} (\vb{x})$.
  We will show $F_1 (\vb{x}) = \beta$.
  First, since every nonempty subset of $\NWC{m}$ has a minimal
  element, there is $\vb{a} \in \NWC{m}$ which is minimal with
  $\ext{F} (\vb{a}) \le \beta$ and $\vb{a} \le \vb{x}$.
  Then since $F_1$ is antitone, we have
  $F_1 (\vb{x}) \le F_1 (\vb{a}) = \ext{F} (\vb{a}) \le \beta$.
  Now suppose that $F_1 (\vb{x}) < \beta$. Let $\gamma := F_1 (\vb{x})$.
  Then
  \[
       \ext{F}(\vb{x}) = \beta \not\le F_1 (\vb{x}) = \gamma.
  \]
  By
  Theorem~\ref{thm:finrep1}, there is a finite subset
  $G \subseteq F$ with $F_{G} = F$. Now by
  Lemma~\ref{lem:fga}\eqref{it:fg2}, the set $A$ consisting of the
  minimal elements of $\{ \vb{c} \in \NWC{m} \mid \ext{F} (\vb{c}) \le \gamma \}$
  satisfies $A \subseteq \NW{m}$. Then we have
  $\{ \vb{c} \in \NWC{m} \mid \ext{F} (\vb{c}) \not\le \gamma \}
   = \NWC{m} \setminus (A{\uparrow})$. We can use
  Proposition~\ref{pro:max} for finding
  a  maximal $\vb{c} \in \NWC{m} \setminus (A{\uparrow})$ with $\vb{c} \ge \vb{x}$.
  Then  $\ext{F} (\vb{c}) \not\le \gamma$.
  Since $\vb{c} \in P$, $F_1 (\vb{c}) = \ext{F} (\vb{c}) \not\le \gamma$.
  Since $F_1$ is antitone, we have $F_1 (\vb{c}) \le F_1 (\vb{x})$
  and thus $F_1 (\vb{x}) \not\le \gamma$, which yields
  $\gamma \not\le \gamma$. This contradiction shows that
  $F_1 (\vb{x}) = \beta$.
\end{proof}
\begin{de}
  Let $m \in \N$, let $\ob{L}$ be a complete lattice,
  and let
  $F : \NW{m} \to \ob{L}$ be an antitone function.
   A \emph{complete representation} of $F$
   is a finite subset $H$ of $\ext{F}$ such that
   $\ext{F}$ is the only antitone function from
   $\NWC{m}$ to $\ob{L}$ with $H \subseteq \ext{F}$.
\end{de}
\begin{thm} \label{thm:hatF}
   Let $m \in \N$, let $\ob{L}$ be a finite lattice,
  and let
  $F : \NW{m} \to \ob{L}$ be an antitone function.
  Then $F$ has a complete representation.
\end{thm}
\begin{proof}
  By assumption,  $\ob{L}$ is finite. For every $\alpha$, the
  set $P_{\alpha}$ defined in Lemma~\ref{lem:unique} is the union
  of two antichains of $\NWC{m}$, and therefore finite
  by Lemma~\ref{lem:dick}.
  Hence the set $P$  defined in Lemma~\ref{lem:unique} is finite,
  and thus by this Lemma $H := \{ (\vb{a}, \ext{F} (\vb{a}) \mid \vb{a} \in P\}$
  satisfies the required property.
\end{proof}

\section{Computing canonical representations} \label{sec:compcan}
In this section, we will see that starting from a finite subset
$G$ of $\NW{m} \times \ob{L}$, we can compute the
canonical representation of $F_G$. When determining further properties
of an antitone functions in Section~\ref{sec:prop},
we will often assume that it is
given by its canonical representation.
We will also see that we can compute a complete representation of $F_G$
from its canonical representation.
We will always assume that a finite lattice $\ob{L}$ is given in a
way that allows to compute with its elements: in particular,
we assume that we can compute a list of the (finitely many) elements
of $\ob{L}$, that we can compute the join and meet of two
elements, and that we can decide whether two elements
of $\ob{L}$ are equal.

For an antitone $F : \NW{m} \to \ob{L}$, and $\alpha \in \ob{L}$,
let
\[
U_F (\alpha) := \{ \vb{a} \in \NW{m} \mid F (\vb{a}) \le \alpha \}.
\]
The next Lemma will allow us to compute $U_{F_G} (\alpha)$ when
$G$ is a finite subset of $\NW{m} \times \ob{L}$.
\begin{lem} \label{lem:ufg}
  Let $\ob{L}$ be a complete lattice, and
  let $G \subseteq \NW{m} \times \ob{L}$.
  Then for each $\alpha \in \ob{L}$, we have
  \begin{equation} \label{eq:ufg}
     U_{F_G} (\alpha) =
     \bigcup \,\,
     \big\{ \bigcap_{ (\vb{b}, \beta) \in J }
           \vb{b}{\uparrow}
            \,\,\, \mid \,\,\,
       J \text{ is a subset of } G \text{ with } (\bigmeet_{(\vb{b}, \beta) \in J} \beta) \le \alpha \big\} .
  \end{equation}
\end{lem}
\begin{proof}
  For $\subseteq$, we choose $\vb{a} \in \NW{m}$ with
  $F_G (\vb{a}) \le \alpha$.
  By the definition of $F_G$, we have
  $F_G (\vb{a}) = \bigmeet_{(\vb{b},\beta) \in G, \,
    \vb{b} \le \vb{a}} \beta$.
  Hence for $J :=
    \{(\vb{b},\beta) \in G, \mid
    \vb{b} \le \vb{a} \}$,
    we have $(\bigmeet_{(\vb{b},\beta) \in J} \beta)  \le \alpha$
    and $\vb{a} \in \bigcap_{(\vb{b},\beta) \in J} \vb{b}{\uparrow}$.

    For $\supseteq$, we choose $\vb{a}$ in the right hand side
    of~\eqref{eq:ufg}. Then there is $J \subseteq G$
    with $(\bigmeet_{(\vb{b}, \beta) \in J} \beta) \le \alpha$ and
    $\vb{a} \in \bigcap_{(b, \beta) \in J} \vb{b}{\uparrow}$.
    This last condition implies that for each $(\vb{b}, \beta) \in J$,
    we have $\vb{b} \le \vb{a}$.
    Hence $J \subseteq \{ (\vb{b}, \beta) \in G \mid
    \vb{b} \le \vb{a}\}$.
    Thus $\alpha \ge \bigmeet_{(\vb{b}, \beta) \in J} \beta
    \ge \bigmeet \{ \beta \mid
    (\vb{b}, \beta) \in G \text{ and }
    \vb{b} \le \vb{a}\} = F_G (\vb{a})$, and therefore
    $\vb{a} \in U_{F_G} (\alpha)$.
\end{proof}
In the sequel, we will compute unions and intersections of upward
closed subsets of $\NW{m}$. Every upward closed subset $U$ of
$\NW{m}$ can be represented as $A{\uparrow}$ for some
finite subset $A$ of $U$: one choice for $A$ is the set
$\Min (U)$ of minimal elements of $U$. Since $\Min (U)$
is an antichain and all antichains in $\NW{m}$ are finite
(Lemma~\ref{lem:dick}), $\Min (U)$ is finite, and since
$\NW{m}$ satisfies the descending chain condition,
every element of $U$ is greater than or equal to some
minimal element of $U$.
When $U \subseteq \NW{m}$ is given as $B{\uparrow}$ for some finite subset
$B$ of $\NW{m}$, then $\Min (B) = \Min (B{\uparrow})$. From
this we obtain that for two finite subsets $A$ and $B$,
$A{\uparrow} = B{\uparrow}$ holds if and only if
$\Min (A) = \Min (B)$. This allows us to decide whether
$A{\uparrow}$ and $B{\uparrow}$ are equal.

Now the union of two upward closed subsets of $\NW{m}$ can
be computed using
\[
    (A{\uparrow}) \cup (B{\uparrow})
    =
    (A \cup B){\uparrow},
\]
and therefore $(A{\uparrow}) \cup (B{\uparrow}) = (\Min (A \cup B)){\uparrow}$.
For computing intersections, we recall that
    for $\vb{a}, \vb{b} \in \NW{m}$, $\vb{a} \sqcup \vb{b}$
    denotes $(\max(a_1, b_1), \ldots, \max(a_m, b_m))$.
    Now for finite subsets $A$ and $B$ of  $\NW{m}$,
    we have
    \[
    (A{\uparrow}) \cap (B{\uparrow}) =
     \{ \vb{a} \sqcup \vb{b} \mid \vb{a} \in A,
     \vb{b} \in B \}{\uparrow},
   \]
   and therefore the set $C$ of  minimal elements of the finite
   set $\{ \vb{a} \sqcup \vb{b} \mid \vb{a} \in A, \vb{b} \in B \}$ satisfies
   $C{\uparrow} =   (A{\uparrow}) \cap (B{\uparrow})$.
   This allows us to compute intersections and unions
   of upward closed sets that are given
   as $A{\uparrow}$ for some finite $A \subseteq \NW{m}$,
   and to return as a result
   a finite subset of $C$ of $\NW{m}$ such that
   $C{\uparrow} =   (A{\uparrow}) \cup (B{\uparrow})$
   (or $C{\uparrow} = (A{\uparrow}) \cap (B{\uparrow})$, resp.).

 \begin{lem} \label{lem:min}
   Let $m \in \N$, let $\ob{L}$ be a finite lattice, and let $\alpha \in \ob{L}$.
   Then there
   is an algorithm whose input is
   a finite subset $G$ of $\NW{m} \times \ob{L}$ and whose output
   is the set
   \(
   \Min\, ( \{ \vb{a} \in \NW{m} \mid F_G (\vb{a}) \le \alpha \} ).
   \)
 \end{lem}
       \begin{proof}
      Lemma~\ref{lem:ufg} represents
      $U_{F_G} (\alpha) = \{ \vb{a} \,\mid\, F_G (\vb{a}) \le \alpha\}$
      as a union of intersections of sets of the form
      $\vb{b}{\uparrow}$ where $\vb{b} \in \NW{m}$.
      By the remarks preceding Lemma~\ref{lem:min}, we
      can therefore compute
    $\Min (U_{F_G} (\alpha)) = \Min (\{ \vb{a} \,\mid\, F_G (\vb{a}) \le \alpha\})$.
       \end{proof}
In a lattice $\ob{L}$, we write $\alpha \lcover \beta$ and say that
\emph{$\beta$ covers $\alpha$} if $\alpha < \beta$ and there is no
$\gamma \in \ob{L}$ with $\alpha < \gamma < \beta$.
    \begin{thm} \label{thm:ccan}
       Let $\ob{L}$ be a finite lattice, and let $G$ be a finite
      subset of $\NW{m} \times \ob{L}$.
      Then there is
      an algorithm to compute the canonical
      representation of $F_G$.
    \end{thm}
    \begin{proof}
      Lemma~\ref{lem:min} allows to compute
          $\Min (U_{F_G} (\alpha)) = \Min (\{ \vb{a} \,\mid\, F_G (\vb{a}) \le \alpha\})$
    for each $\alpha \in \ob{L}$.
    In order to be able  compute the canonical representation of $F_G$,
    we show that for every antitone function $F$, we have
    \begin{multline}\label{eq:eqandleq}
      \Min (\{ \vb{a} \,\mid\, F (\vb{a}) = \alpha\}) \\
      =
      \Min (\{ \vb{a} \,\mid\, F (\vb{a}) \le \alpha\}) \setminus
      ( \bigcup_{\beta \in \ob{L}, \beta \lcover \alpha} \Min (\{ \vb{a} \,\mid\,
      F (\vb{a}) \le \beta \} )).
    \end{multline}
    For $\subseteq$, let $\vb{a}$ be minimal with
    $F( \vb{a} ) = \alpha$.
    To show that $\vb{a}$ is minimal with $F(\vb{a}) \le \alpha$,
    we choose $\vb{b} \le \vb{a}$ with
    $F(\vb{b}) \le \alpha$.  Then by antitony $F(\vb{b}) \ge F(\vb{a})
    = \alpha$,
    and therefore $F(\vb{b}) = \alpha$. The minimality of
    $\vb{a}$ now yields $\vb{b} = \vb{a}$.
    For all $\beta \lcover \alpha$, we have
    $F(\vb{a}) = \alpha \not\le \beta$, and therefore
    $\vb{a}$ lies in the right hand side of~\eqref{eq:eqandleq}.

    For $\supseteq$, we choose $\vb{a}$ in the right hand side
    of~\eqref{eq:eqandleq}. We first show $F(\vb{a}) = \alpha$.
    Seeking a contradiction, we suppose $F(\vb{a}) < \alpha$.
    Then there is $\beta \lcover \alpha$ with
    $F (\vb{a}) \le \beta$. We now show that $\vb{a}$ is
    minimal with $F(\vb{a}) \le \beta$. To this end,
    let $\vb{b} \le \vb{a}$ with $F(\vb{b}) \le \beta$.
    Then $F(\vb{b}) \le \alpha$, and since $\vb{a}$ is minimal
    with $F(\vb{a}) \le \alpha$, we have $\vb{b} = \vb{a}$.
    Thus $\vb{a}$ is
    minimal with $F(\vb{a}) \le \beta$.
    Therefore, $\vb{a}$ is not an element of the right
    hand side; this contradiction proves $F (\vb{a}) = \alpha$.
    We now prove that $\vb{a}$ is minimal with
    $F(\vb{a}) = \alpha$; this follows from the fact that
    $\vb{a}$ is even minimal with $F(\vb{a}) \le \alpha$.
    This completes the proof of~\eqref{eq:eqandleq}

    We use~\eqref{eq:eqandleq} to compute, for each $\alpha \in \ob{L}$, the set
    $\Min (\{ \vb{a} \,\mid\, F (\vb{a}) = \alpha\})$.
    For computing the right hand side of~\eqref{eq:eqandleq},
    Lemma~\ref{lem:min} allows us to compute
    $\Min (\{ \vb{a} \,\mid\, F_G (\vb{a}) \le \gamma\})$
    for each $\gamma \in \ob{L}$. Since we can
    compute in $\ob{L}$, we can find those
    $\beta \in \ob{L}$ with
    $\beta \lcover \alpha$. What remains to obtain the
    right hand side of~\eqref{eq:eqandleq} is to perform
    set theoretic operations with finite sets.

    Now $G := \{ (\vb{a}, \alpha) \mid \alpha \in \ob{L}, \,
    \vb{a} \text{ is a minimal element of } F^{-1}(\{\alpha\}) \}$
    is the required canonical representation of $F$.
\end{proof}

    \section{Computing complete representations}
In this section, we compute a complete representation of
$F_G$. Here, given $G$ our goal is to find a finite
subset $H$ of $\NWC{m} \times \ob{L}$ such that
$\ext{F_G}$ is the only antitone function with
$H \subseteq \ext{F_G}$.
\begin{lem} \label{lem:max2}
  Let $A$ be a finite subset of $\NW{m}$, and let
  $\vb{b}$ be a maximal element of $\NWC{m} \setminus (A{\uparrow})$,
  and let $i \in \ul{m}$. Then $b_i \in \{ a_i - 1 \mid \vb{a} \in A, a_i > 0 \}
  \cup \{ \infty \}$.
\end{lem}
\begin{proof}
  Suppose $b_i \neq \infty$. By the maximality of $\vb{b}$,
  $\vb{b} + \vb{e}_i  \in A{\uparrow}$, and therefore
  there is $\vb{a} \in A$ with
  $\vb{a} \le \vb{b} + \vb{e}_i$.
  Hence $a_i \le b_i + 1$. If $a_i \le b_i$, then
  using $\vb{a} \le \vb{b} + \vb{e}_i$ we obtain that
  $\vb{a} \le \vb{b}$, which implies $\vb{b} \in A{\uparrow}$,
  a contradiction.
  Thus $a_i = b_i + 1$, and therefore $a_i > 0$ and $b_i = a_i - 1$.
\end{proof}
\begin{lem} \label{lem:max}
   Let $m \in \N$, let $\ob{L}$ be a finite lattice, and let $\alpha \in \ob{L}$.
   Then there
   is an algorithm whose input is
   a finite subset $G$ of $\NW{m} \times \ob{L}$ and whose output
   is the set
   \(
   \Max\, ( \{ \vb{a} \in \NWC{m} \mid \ext{F_G} (\vb{a}) \not\le \alpha \} ).
   \)
 \end{lem}
\begin{proof}
  We first use Lemma~\ref{lem:min} to compute the set
  \[
   A  := \Min\, ( \{ \vb{a} \in \NW{m} \mid F_G (\vb{a}) \le \alpha \} ).
  \]
  By Lemma~\ref{lem:fga}\eqref{it:fg2},
  $\{ \vb{b} \in \NWC{m} \mid \ext{F_G} (\vb{b}) \le \alpha \}
  =
   A{\uparrow} = \{ \vb{b} \in \NWC{m} \mid \exists \vb{a} \in A : \vb{a} \le \vb{b} \}$.
   Using Lemma~\ref{lem:max2}, we find a finite subset
   $B$ of $\NWC{m}$ that contains all maximal elements of
   $\NWC{m} \setminus (A{\uparrow})$. Now the set
   $\Max (B)$ is the set of maximal elements of
   $\{ \vb{b} \in \NWC{m} \mid \ext{F_G} (\vb{b}) \not\le \alpha \}$.
 \end{proof}
\begin{thm} \label{thm:ccomplete}
       Let $\ob{L}$ be a finite lattice, and let $G$ be a finite
      subset of $\NW{m} \times \ob{L}$.
      Then there is
      an algorithm to compute a complete representation
      of $F_G$.
\end{thm}
\begin{proof}
  Proceeding as in the proof of Theorem~\ref{thm:ccan},
  we fix
  $\alpha \in \ob{L}$ and use Lemma~\ref{lem:ufg} to
  compute the set $A_{\alpha}$
  of minimal elements of
  $\{ \vb{a} \in \NW{m} \mid F_G (\vb{a}) \le \alpha \}$.
  By Lemma~\ref{lem:fga}\eqref{it:fg2},
  $\{ \vb{b} \in \NWC{m} \mid \ext{F_G} (\vb{b}) \le \alpha \}
  =
   A_{\alpha}{\uparrow}$.
   Using Lemma~\ref{lem:max2}, we can compute a finite subset
   $B_\alpha$ of $\NWC{m}$ that contains all maximal elements of
   $\NWC{m} \setminus (A{\uparrow})$.
   We set $P := \bigcup_{\alpha \in \ob{L}} (A_{\alpha} \cup B_{\alpha})$
   and $P' := \{ (\vb{a}, F_G(\vb{a}) \mid \vb{a} \in P \}$;
   the evaluation of $\ext{F_G}$ at $\vb{a}$ can be accomplished
   using Lemma~\ref{lem:fga}\eqref{it:fg1}.
       From Lemma~\ref{lem:unique}, we obtain that
   $\ext{F}$ is the only antitone function
   with $P' \subseteq \ext{F}$.
\end{proof}

Next, we seek to determine whether there is exactly one
antitone function through a given set of points.
\begin{lem} \label{lem:decidecomplete}
      Let $\ob{L}$ be a complete lattice, let $G$ be a finite
      subset of $\NW{m} \times \ob{L}$, and let
      $H \subseteq \ext{F_G}$. Then the following are equivalent:
      \begin{enumerate}
      \item \label{it:com1}
   $\ext{F_G}$ is the only antitone function from
   $\NWC{m}$ to $\ob{L}$ that contains $H$ as a subset.
      \item \label{it:com2}
        For all $\alpha \in \ob{L}$ and
        for all $\vb{b} \in \Min (\{ \vb{a} \in \NW{m} \mid F_G (\vb{a}) \le \alpha \})$, we have
         \begin{equation} \label{eq:meet}
        \bigmeet \{ \delta \mid (\vb{d}, \delta) \in H \text{ and } \vb{d} \le \vb{b} \} \le \alpha,
        \end{equation}
        and for all $\alpha \in \ob{L}$ and
        for all $\vb{b} \in \Max (\{ \vb{a} \in \NWC{m} \mid \ext{F_G} (\vb{a}) \not\le \alpha \})$,
        we have
         \begin{equation} \label{eq:vee}
         \bigvee \{ \delta \mid (\vb{d}, \delta) \in H \text{ and } \vb{b} \le \vb{d} \} \not\le \alpha.
        \end{equation}
      \end{enumerate}
\end{lem}
\begin{proof}
  \eqref{it:com2}$\Rightarrow$\eqref{it:com1}:
  Let $F$ be an antitone function from $\NWC{m}$ to $\ob{L}$
  with $H \subseteq F$. We show $F = \ext{F_G}$. To this
  end, we fix $\vb{x} \in \NWC{m}$.
  Let $\alpha := \ext{F_G} (\vb{x})$. We choose $\vb{b}$ minimal
  with $\vb{b} \le \vb{x}$ and $\ext{F_G} (\vb{b}) \le \alpha$.
  Then by Lemma~\ref{lem:fga}\eqref{it:fg2},  $\vb{b} \in \NW{m}$.
  Since $F$ is antitone, we know that $F (\vb{b}) \le F (\vb{d})$
  for all $\vb{d} \in \NW{m}$ with $\vb{d} \le \vb{b}$.
  Since $H \subseteq F$, we therefore have
  \[
    F (\vb{b}) \le
  \bigmeet \{ \delta \mid (\vb{d}, \delta) \in H \text{ and } \vb{d} \le \vb{b} \} \le \alpha.
  \]
  Since $F$ is antitone, we have $F(\vb{x}) \le F (\vb{b})$,
  and thus $F(\vb{x}) \le \alpha$.
  Now assume that $F(\vb{x}) < \alpha$.
  Then for $\beta := F (\vb{x})$, we have
  $\ext{F_G} (\vb{x}) \not\le \beta$.
  The set $A$ that contains the
  minimal elements of $\{ \vb{a} \in \NWC{m} \mid \ext{F_G} (\vb{a}) \le \beta \}$
  satisfies $A \subseteq \NW{m}$. Then we have
  $\{ \vb{a} \in \NWC{m} \mid \ext{F_G} (\vb{a}) \not\le \beta \}
  = \NWC{m} \setminus (A{\uparrow})$.
  By
  Proposition~\ref{pro:max}, there is
  a  maximal $\vb{b} \in \NWC{m} \setminus (A{\uparrow})$ with $\vb{b} \ge \vb{x}$.
  Then  $\ext{F_G} (\vb{b}) \not\le \beta$, and therefore by assumption
  $\bigvee \{ \delta \mid (\vb{d}, \delta) \in H, \vb{b} \le \vb{d} \}
   \not\le \beta$.

  Since $F$ is antitone, we know that $F (\vb{b}) \ge F (\vb{d})$
  for all $\vb{d} \in \NWC{m}$ with $\vb{b} \le \vb{d}$.
  Since $H \subseteq F$, we therefore have
  \[
    F (\vb{b}) \ge
  \bigjoin \{ \delta \mid (\vb{d}, \delta) \in H \text{ and } \vb{b} \le \vb{d} \}.
  \]
  Thus by the assumption, $F(\vb{b}) \not\le \beta$.
  Since $F$ is antitone, we have  $F(\vb{x}) \ge F (\vb{b})$,
  and therefore $F (\vb{x}) \not\le \beta$. This yields
  $\beta \not\le \beta$, a contradiction. Thus
  $F (\vb{x}) = \alpha$.

   \eqref{it:com1}$\Rightarrow$\eqref{it:com2}:
   We define two antitone functions $F_1$ and $F_2$ by
   \[
   \begin{array}{rcl}
        F_1 (\vb{x}) & := &
        \bigmeet \{ \delta \mid (\vb{d}, \delta) \in H \text{ and } \vb{d} \le \vb{x} \}, \\
        F_2 (\vb{x}) & := &  \bigvee \{ \delta \mid (\vb{d}, \delta) \in H \text{ and } \vb{x} \le \vb{d} \}
   \end{array}
   \]
   for all
   $\vb{x} \in \NWC{m}$.
   We first show $
   H \subseteq F_1.
   $
   To this end, let $(\vb{b}, \beta) \in H$.
   Then
   $F_1 (\vb{b}) = \bigmeet D$,  where $D = \{ \delta \mid (\vb{d}, \delta) \in H \text{ and } \vb{d} \le \vb{b} \}$. Then $\beta \in D$, and
   since $H \subseteq \ext{F_G}$ and $\ext{F_G}$ is antitone,
   we have $\beta \le \delta$ for all $(\vb{d}, \delta) \in H$ with
   $\vb{d} \le \vb{b}$, hence $\beta \le \delta$ for all $\delta \in D$.
   Therefore $\bigmeet D = \beta$, and thus $F_1 (\vb{b}) = \beta$.
   A dual argument shows $H \subseteq F_2$.
   Therefore, $\ext{F_G} = F_1 = F_2$.

   Now let $\alpha \in \ob{L}$, and let $\vb{b} \in \NW{m}$
   be minimal with $\ext{F_G} (\vb{b}) \le \alpha$.
   Then $F_1 (\vb{b}) \le \alpha$, which implies~\eqref{eq:meet}.
   For proving~\eqref{eq:vee}, we let
   $\vb{b} \in \NWC{m}$ be maximal with $\ext{F_G} (\vb{b}) \not\le \alpha$.
   Then $F_2 (\vb{b}) \not\le \alpha$, which implies~\eqref{eq:vee}.
\end{proof}

\begin{thm} \label{thm:deccomplete}
  Let $m \in \N$, let $\ob{L}$ be a finite lattice. Then there
  is an algorithm that, given the finite sets
  $G \subseteq \NW{m} \times \ob{L}$ and
  $H \subseteq \NWC{m} \times \ob{L}$, decides whether
  $H$ is a complete representation of $\ext{F_G}$.
\end{thm}
\begin{proof}
  Using Lemma~\ref{lem:fga}\eqref{it:fg1}, we can evaluate
  $\ext{F_G}$, and in this way we can check whether $\ext{F_G}$
  passes through all points of $H$, i.e., whether
  $H \subseteq \ext{F_G}$. If this condition is fulfilled,
  we check whether condition~\eqref{it:com2} of
  Lemma~\ref{lem:decidecomplete} is satisfied. To this end,
  we use Lemmas~\ref{lem:min}~and~\ref{lem:max} to compute,
  for each $\alpha \in \ob{L}$, the finite sets
  $\Min (\{ \vb{a} \in \NW{m} \mid F_G (\vb{a}) \le \alpha \})$
  and $\Max (\{ \vb{a} \in \NWC{m} \mid \ext{F_G} (\vb{a}) \not\le \alpha \})$.
  With these sets available, we can check whether~\eqref{eq:meet}~and~\eqref{eq:vee} are satisfied, and exactly in this case, we return ``true''.
  \end{proof}

We conclude this section with an example of these
representations. In this example,
the lattice $\ob{L} = \{1,2,4,13,26,52\}$
is the set of divisors of $52$, and $x \meet y := \gcd(x,y)$,
and we choose \[G = \{ ((10, 20), 26), ((30, 5), 4)) \}. \]
  Then
  the canonical representation of $F_G : \NW{2} \to \ob{L}$ is
  \[G_1 = \{ ((0, 0), 52), ((10, 20), 26), ((30, 5), 4), ((30, 20), 2) \}, \]
  and a complete representation of $F_G$ is
  \begin{multline*}
     G_2 =\{ ((0, 0), 52), ((10, 20), 26), ((30, 5), 4), ((30, 20),
  2), ((9, \infty), 52), \\ ((29, 19), 52), ((29, \infty),
  26), ((\infty, 4), 52), ((\infty, 19),
  4), ((\infty, \infty), 2) \}.
  \end{multline*}

 \section{Learning antitone functions}

 In the next section, we will see that if we are given a canonical
 representation of the antitone function that encodes an operation
 sequence with \HC{3} and \HC{4}, we can decide from this canonical
 representation whether the operation sequence satisfies
 \HC{1}, \HC{2}, \HC{7}, or (\HC{2} and \HC8).
 This raises the question how to find a canonical, or even some,
 representation, for a given antitone function $F$.
 Here we suppose that $F$ is given as a ``black box'': we can
 evaluate $F$ at arbitrary places $\vb{a} \in \NW{m}$, but
 we should come to a conclusion after finitely many such
 evaluations. Hence the goal would be to evaluate $F$ at enough points
 to allow us to  determine a finite set $G$ with $F = F_G$.
 However, from
 Proposition~\ref{pro:nonunique}, we gather that our function
 $F$ may be such that for every finite subset $H$ of $F$,
 there is a second antitone function $F' \not= F$
 with $H \subseteq F'$.
 Hence evaluating $F$ at a finite number of points may
 not result in a subset $G$ with $F = F_G$.
 However, if we are allowed to evaluate the extension $\ext{F}$ also at points
 having $\infty$ in one of their components, then there is
 an algorithm to find $G \subseteq \NW{m} \times \ob{L}$ with $F = F_G$.
 We prepare for this algorithm with two lemmas.
 \begin{lem} \label{lem:hatfb}
   Suppose $F_1, F_2$ are antitone functions, and suppose
   that $\vb{a} \in \NWC{m}$ is such that
   $\ext{F_1} (\vb{a}) < \ext{F_2} (\vb{a})$.
   Then there exists
   $\vb{b} \in \NW{m}$ with
   $\ext{F_1} (\vb{b}) < \ext{F_2} (\vb{b})$.
 \end{lem}
 \begin{proof}
   By Lemma~\ref{lem:fga}\eqref{it:fg15}, there is $\vb{b} \in \NW{m}$
   with $\vb{b} \le \vb{a}$ and
 $F_1 (\vb{b}) = \ext{F_1} (\vb{a})$.
    Then $\ext{F_1} (\vb{b}) = F_1 (\vb{b}) = \ext{F_1} (\vb{a}) < \ext{F_2} (\vb{a})$,
    which, by the antitony of $\ext{F_2}$, is $\le \ext{F_2} (\vb{b})$.
 \end{proof}
When $F_1$ and $F_2$ are given by ``black boxes'' that
 yield the value of $F_1 (\vb{a})$ and $F_2 (\vb{a})$ for each
 input $\vb{a} \in \NW{m}$,
 and we know that $F_1 (\vb{a}) < F_2 (\vb{a})$ for some
 $\vb{a} \in \NWC{m}$, then $\vb{b} \in \NW{m}$
 with $F_1 (\vb{b}) < F_2 (\vb{b})$ can be found algorithmically:
 we enumerate all $\vb{b} \in \NW{m}$
 and stop when a $\vb{b}$ with $F_1 (\vb{b}) < F_2 (\vb{b})$ is found.

 The second lemma was essentially already proved in \cite{AM:SOCO}.
 \begin{lem}[cf. {\cite[Theorem~1.2]{AM:SOCO}}]
   Let $\ob{L}$ be a finite lattice, and let $m \in \N$.
   We order antitone functions from $\NW{m}$ to $\ob{L}$ by $F_1 \le F_2$ if
   $F_1 (\vb{a}) \le F_2 (\vb{a})$ for all $\vb{a} \in \NW{m}$.
   Then there is no infinite descending chain of antitone functions.
 \end{lem}
 \begin{proof}
   For each antitone function, the set
   $A (F) := \{ (\vb{a}, \alpha) \mid F (\vb{a}) \le \alpha \}$ is an
   upward closed subset of $\NW{m} \times \ob{L}$.
   Suppose that $F_1 > F_2 > F_3 > \ldots$ is an infinite descending
   chain of antitone functions, then
   $A(F_1) \subseteq A (F_2) \subseteq A(F_3) \subseteq \ldots$
   is an infinite ascending chain of upward closed subsets of
   $\NW{m} \times \ob{L}$. By Lemma~\ref{lem:dick},
   $A := \bigcup_{i \in \N} A (F_i)$ has only finitely many
   minimal elements, which yields $j \in \N$ with
   $A \subseteq A(F_j)$. Hence $A (F_{j+1}) = A (F_j)$.
   Since $(F_i)_{i \in \N}$ is descending, there is $\vb{a} \in \NW{m}$ such that
   $F_{j+1} (\vb{a}) < F_{j} (\vb{a})$.
   Then $(\vb{a}, F_{j+1} (\vb{a})) \in A(F_{j+1})$ and
   $(\vb{a}, F_{j+1} (\vb{a}))  \not\in A(F_{j})$, contradicting
   $A (F_{j+1}) = A (F_j)$.
   Hence there is no infinite descending chain of antitone functions.
 \end{proof}
 The following Theorem tells us how to ``learn'' an antitone
 function $F$  (in the sense of computational learning theory)
 from finitely many evaluations of its extension $\ext{F}$.
 \begin{thm} \label{thm:learn}
   Let $\ob{L}$ be a finite lattice,
   let $m \in \N$,
   and let $F : \NW{m} \to \ob{L}$
   be an antitone function.
   Suppose that we can compute the value $\ext{F} (\vb{a})$ for each
   $\vb{a} \in \NWC{m}$.
   Then there is an algorithm that computes a finite subset
   $G$ of $\NW{m} \times \ob{L}$ such that $F = F_G$.
 \end{thm}
 \begin{proof}
   We consider the following procedure:
      \begin{algorithmic}[1]
        \State $G \gets \emptyset$
        \State \text{compute a complete representation $H$ of $F_G$}
        \While{$H \not\subseteq \ext{F}$}
                \State \text{find $\vb{a} \in \NW{m}$ with
                  $F (\vb{a}) < F_G (\vb{a})$}
                  \State $G \gets G \cup \{ (\vb{a}, F (\vb{a})) \}$
                  \State  \text{compute a complete representation $H$ of $F_G$}
        \EndWhile
              \State \text{return $G$.}
      \end{algorithmic}
      We first explain how each of these steps can be computed:
      A complete representation $H$ of $F_G$ can be computed
      using Theorem~\ref{thm:ccomplete}.
      Throughout the algorithm, the property
      $G \subseteq F$ is invariant. Hence by Lemma~\ref{lem:largest},
      we always have $\forall \vb{a} \in \NW{m} : F(\vb{a}) \le F_G (\vb{a})$
      and $H \subseteq \ext{F_G}$.
      The condition $H \subseteq \ext{F}$ can be checked by
      evaluating $\ext{F}$ at finitely many points
      $\{ \vb{a} \mid (\vb{a}, \alpha) \in H \}$, which we
      can do by assumption.
      If $H \not\subseteq \ext{F}$, then $\ext{F} \neq \ext{F_G}$.
      Since $F \le F_G$, we then have
      $\ext{F} < \ext{F_G}$, and therefore by Lemma~\ref{lem:hatfb}
      and the remarks following this Lemma, we can find
      an
      $\vb{a} \in \NW{m}$ with
      $F (\vb{a}) < F_G (\vb{a})$.
      For evaluating $F_G (\vb{a})$,
      we use Definition~\ref{de:fga}; here we use the
      assumption that $\ob{L}$ is ``computable'' and thus we are
      able
      to compute the meet that appears in~\eqref{eq:fga}.

        This search for $\vb{a}$ can
      be avoided if $H \cap (\NW{m} \times \ob{L}) \not\subseteq \ext{F}$:
      then every $\vb{a} \in (H \cap (\NW{m} \times \ob{L})) \setminus \ext{F}$ satisfies
      $\ext{F_G} (\vb{a}) \neq F(\vb{a})$ and thus
      $F(\vb{a}) < F_G (\vb{a})$.

      Next, we will show that the algorithm terminates.
      The reason is that $F_G$ drops in every execution of the
      body of the while loop.
      If $G_1 = G$ and $G_2 = G_1 \cup \{ (\vb{a}, F(\vb{a})) \}$ and
      $F (\vb{a}) \le F_{G_1} (\vb{a})$, then
      $F_{G_2} \le F_{G_1}$.
      This holds because $G_1 \subseteq G_2 \subseteq F_{G_2}$,
      and $F_{G_1}$ is the largest function containing $G_1$.
      Furthermore, we have $F_{G_2} (\vb{a}) = F (\vb{a})$
      (since $G_2 \subseteq F$, this is guaranteed by
      Lemma~\ref{lem:largest}\eqref{it:l3}) and
      $F (\vb{a}) < F_G (\vb{a}) = F_{G_1} (\vb{a})$.
       Now since there is no infinite descending chain
      of antitone functions, the while loop must terminate,
      which implies that then the condition $H \subseteq \ext{F}$
      is satisfied.

      We will now show that the algorithm produces a correct result.
      Since when the algorithm ends, $H \subseteq \ext{F}$,
      and $H \subseteq \ext{F_G}$ holds throughout the algorithm,
      the fact that $H$ is a complete representation yields
      $\ext{F} = \ext{F_G}$.
 \end{proof}
 We notice that if we return $H$ instead of $G$ in the last
 line of our procedure, then the procedure computes a complete
 representation of $F$. Let us give one example to show that this algorithm may
 be of  considerable complexity if we make no assumptions
 on the input $F$: we take $m=1$ and  $F : \N_0 \to \{0,1\}$ with
 $F (x) = 1$ for $x \le 10^{10^{10}}$, and
 $F (x) = 0$ for $x >10^{10^{10}}$. Then in the first
 execution of the while loop,
 we have to find an $a \in \N$ with $F (a) < F_G(a) = 1$,
 which we find only when we evaluate $F$ at values
 exceeding $10^{10^{10}}$.

\section{Deciding properties of an antitone function} \label{sec:prop}
In Lemma~\ref{lem:HCencoding}, several conditions on the commutator
sequence were translated into conditions on their encoding. In this
section, we give versions that can effectively be tested when
we are given a canonical representation of an antitone function.
The main point is to reduce checking a property for all $\vb{a} \in \NW{m}$
to checking it at some finite set of places.
\begin{thm}
  Let $\ob{L}$ be a lattice, and let
  $F : \NW{m} \to \ob{L} = \{\lambda_1, \ldots, \lambda_m \}$ be an antitone function,
  and let $G$ be the canonical representation of $F$.
  Then we have
  \begin{enumerate}
  \item \label{it:hc1} $F$ satisfies (HC1) if and only if
    for all $j \in \ul{m}$, we have
    $F (\vb{e}_j) \le \lambda_j$.
  \item \label{it:hc2} $F$ satisfies (HC2) if and only if
    for all $i,j \in \ul{m}$ with
    $\lambda_i \le \lambda_j$ and for all
    $(\vb{b}, \beta) \in G$ with $b_j > 0$,
    we have
    $F (\vb{b} - \vb{e}_j + \vb{e}_i) \le F(\vb{b})$.
  \item \label{it:hc8}
    Suppose that $F$ satisfies (HC2). Then
    $F$ satisfies (HC8) if and only if for all $(\vb{a}, \alpha) \in G$
    and for all $\vb{b} \in \NW{m}$ with $\vb{b} \le \vb{a}$ with
    $F (\vb{b}) = \lambda_j$, we have
    $F(\vb{a} - \vb{b} + \vb{e}_j) \le F (\vb{a})$.
  \end{enumerate}
\end{thm}

\emph{Proof:}
\eqref{it:hc1} Since $F$ is antitone, $F(\vb{e}_j) \le \lambda_j$ implies
$F (\vb{a}) \le \lambda_j$ for all $\vb{a}$ with $a_j \ge 1$.

\eqref{it:hc2} Let $\vb{a} \in \NW{m}$ with $a_j > 0$, and let
$\alpha := F (\vb{a})$. Since $G$ is the canonical representation of
$F$, there is
$\vb{b} \in \NW{m}$ with $(\vb{b}, \alpha) \in G$ and
$\vb{b} \le \vb{a}$.
If $b_j = 0$, then $\vb{b} \le \vb{a} - \vb{e}_j + \vb{e}_i$ and therefore
$F(\vb{a} - \vb{e}_j + \vb{e}_i) \le F (\vb{b}) = \alpha$.
If $b_j > 0$, then
$\vb{b} - \vb{e}_j + \vb{e}_i  \le \vb{a}- \vb{e}_j + \vb{e}_i$, and therefore
$F(\vb{a} - \vb{e}_j + \vb{e}_i) \le
F(\vb{b} - \vb{e}_j + \vb{e}_i) \le F (\vb{b}) = \alpha$.

\eqref{it:hc8}
 Let $j \in \ul{m}$, let
       $\vb{a} \in \NW{m}$, and let $\vb{b} \in \NW{m}$
 with $F(\vb{b}) = \lambda_j$ and $\vb{b} \le \vb{a}$.
 We want to show
 \begin{equation} \label{eq:abe}
 F (\vb{a} - \vb{b} + \vb{e}_j) \le F (\vb{a}).
 \end{equation}
 Let $\alpha := F(\vb{a})$, and
 let $\vb{c}$ be such that $(\vb{c}, \alpha) \in G$ and
 $\vb{c} \le \vb{a}$. Then $\vb{c} \sqcap \vb{b} \, \le \, \vb{b}$,
 and therefore $F (\vb{c} \sqcap \vb{b})  \ge  F (\vb{b}) = \lambda_j$.
 Let $i \in \ul{m}$ be such that $\lambda_i = F (\vb{c} \sqcap \vb{b})$.
 Then by the assumption, $\alpha = F (\vb{c}) \ge F (\vb{c} - (\vb{c} \sqcap \vb{b}) + \vb{e}_i)$.
 Next, we prove
 \begin{equation} \label{eq:ca}
   \vb{c} - (\vb{c} \sqcap \vb{b}) \le \vb{a} - \vb{b}.
 \end{equation}
 We consider the $j$\,th entry.
 If $c_j \le b_j$, then the $j$\,th entry of the left hand side
 is $0$.
 If $b_j \le c_j$, then the $j$\,th entry of the left hand side
 is $c_j - b_j$. Since $c_j \le a_j$, this is at most the $j$\,th
 entry $a_j- b_j$ of the right hand side.
 From~\eqref{eq:ca} and the fact that $F$ is antitone,
 we obtain
 $F (\vb{c} - (\vb{c} \sqcap \vb{b}) + \vb{e}_i) \ge
 F (\vb{a} - \vb{b} + \vb{e}_i)$. Now by (HC2), we obtain
 $F (\vb{a} - \vb{b} + \vb{e}_i) \ge
 F (\vb{a} - \vb{b} + \vb{e}_j)$. Altogether,
 we obtain~\eqref{eq:abe}, which completes the
 proof of~\eqref{it:hc8}. \qed

 Finally, we give a way to decide the join distributivity \HC{7}
 of an antitone function $F_G$.
 \begin{lem} \label{lem:compare}
   Let $G$ be a subset of $\NW{m} \times \ob{L}$, and let
   $F$ be an antitone function.
   Then $F(\vb{a}) \le F_G (\vb{a})$ for all $\vb{a} \in \NW{m}$
   if and only if $F(\vb{b}) \le F_G (\vb{b})$ for all
   $(\vb{b}, \beta) \in G$.
 \end{lem}
 \begin{proof}
   We only have to prove the ``if''-direction. By the assumption,
   $F (\vb{b}) \le F_G (\vb{b}) \le \beta$ for all $(\vb{b}, \beta) \in G$.
   By Lemma~\ref{lem:largest}\eqref{it:l2}, $F_G$ is the largest
   function with $F_G (\vb{b}) \le \beta$ for all $(\vb{b}, \beta) \in G$,
   and thus $F \le F_G$.
 \end{proof}

 The next lemma tells that shifting an antitone function
 commutes with forming its extension. We extend
 $+$ and $-$ by defining $\infty + a = \infty - a = \infty$ for
 all $a \in \N_0$.
  \begin{lem} \label{lem:shift}
    Let $m \in \N$, and let $\ob{L}$ be a finite lattice,
    let $\vb{a} \in \NW{m}$, and let $F : \NW{m} \to \ob{L}$
    be an antitone function.
    We define an antitone function $F_{\vb{a}}$ by
    $F_{\vb{a}} (\vb{x}) := F (\vb{x} + \vb{a})$ for $\vb{x} \in \NW{m}$.
    Then for every $\vb{x} \in \NWC{m}$, we have
    $\ext{F_{\vb{a}}} (\vb{x}) = \ext{F} (\vb{x} + \vb{a})$.
  \end{lem}
  \begin{proof}
    We fix $\vb{x} \in \NWC{m}$ and compute
    $\ext{F_{\vb{a}}} (\vb{x}) =
    \bigmeet \{ F_{\vb{a}} (\vb{y}) \mid \vb{y} \in \NW{m}, \vb{y} \le \vb{x} \}
    =
    \bigmeet \{ F (\vb{y} + \vb{a}) \mid \vb{y} \in \NW{m}, \vb{y} \le \vb{x} \} =: \bigmeet B$ and
    $\ext{F} (\vb{x} + \vb{a}) =
    \bigmeet \{ F (\vb{z}) \mid \vb{z} \in \NW{m}, \vb{z} \le \vb{x} + \vb{a} \} =: \bigmeet C$.
    In order to show $\bigmeet B = \bigmeet C$, we show that
    for all $\beta \in B$, there is $\gamma \in C$ with $\beta \ge \gamma$
    (this establishes $\bigmeet B \ge \bigmeet C$) and that
    for all $\gamma \in C$, there is $\beta \in B$ with $\gamma \ge \beta$
    (this establishes $\bigmeet C \ge \bigmeet B$).
    To this end let $\beta = F (\vb{y} + \vb{a})$ with  $\vb{y} \in \NW{m}$
    and $\vb{y} \le \vb{x}$. Then $\vb{y} + \vb{a} \le \vb{x} + \vb{a}$.
    Hence $\gamma := F(\vb{y} + \vb{a}) \in C$.
    For the other direction, let
    $\gamma =  F (\vb{z})$ with  $\vb{z} \in \NW{m}$ and
    $\vb{z} \le \vb{x} + \vb{a}$. We show that then
    \begin{equation} \label{eq:za}
      \vb{z} - (\vb{z} \sqcap \vb{a}) \le \vb{x}.
    \end{equation}
    We fix $i \in \ul{m}$. If $z_i \le a_i$, then
    $(\vb{z} \sqcap \vb{a})_i = z_i$ and thus
    the $i$\,th component of the left hand side of~\eqref{eq:za}
    is equal to $(\vb{z} - (\vb{z} \sqcap \vb{a}))_i =
    z_i - z_i = 0$, which is $\le x_i$.
    If $z_i > a_i$, then
    $(\vb{z} - (\vb{z} \sqcap \vb{a}))_i =
    z_i - a_i \le (x_i + a_i) - a_i \le x_i$.
    This establishes~\eqref{eq:za}.
    Now we have
    $\beta := F (\vb{z} - (\vb{z} \sqcap \vb{a}) + \vb{a}) \in B$.
    Since $\vb{z} - (\vb{z} \sqcap \vb{a}) + \vb{a} \ge \vb{z}$,
    we have
    $\beta = F (\vb{z} - (\vb{z} \sqcap \vb{a}) + \vb{a})
    \le  F (\vb{z}) = \gamma$.
    Thus $\bigmeet B = \bigmeet C$.
  \end{proof}
  \begin{lem} \label{lem:join}
    Let $m \in \N$, and let $\ob{L}$ be a finite lattice,
    and let $F_1, F_2 : \NW{m} \to \ob{L}$
    be  antitone functions.
    We define an antitone function $F_3$ by
    $F_3 (\vb{x}) := F_1 (\vb{x}) \join F_2 (\vb{x})$
    for $\vb{x} \in \NW{m}$.
    Then for every $\vb{x} \in \NWC{m}$, we have
    $\ext{F_3} (\vb{x}) = \ext{F_1} (\vb{x}) \join \ext{F_2} (\vb{x})$.
  \end{lem}
  \begin{proof}
    We fix $\vb{x} \in \NWC{m}$.
    Since $F_1, F_2, F_3$  are antitone, Theorem~\ref{thm:finrep1}
    yields, for each $i \in \{1,2,3\}$,  a finite set $G_i \subseteq F_i$
    such that $F_i = F_{G_i}$; this observation
    allows us to use Lemma~\ref{lem:fga}\eqref{it:fg15}
    to find $\vb{a}, \vb{b}, \vb{c}  \in \NW{m}$
    with $\vb{a} \le \vb{x}$, $\vb{b} \le \vb{x}$, $\vb{c} \le \vb{x}$
    and
    $F_1 (\vb{a}) = \ext{F_1} (\vb{x})$,
    $F_2 (\vb{b}) = \ext{F_2} (\vb{x})$, $F_3 (\vb{c}) = \ext{F_3} (\vb{x})$.
    Let $\vb{z} := \vb{a} \sqcup \vb{b} \sqcup \vb{c}$.
    Since $\ext{F_1}, \ext{F_2}, \ext{F_3}$ are antitone,
    $\ext{F_i} (\vb{x}) = \ext{F_i} (\vb{z})$ for
    each $i \in \{1,2,3\}$.
    Thus we obtain $\ext{F_3} (\vb{x}) =
    F_3 (\vb{z}) =
    F_1 (\vb{z}) \join F_2 (\vb{z})
    =
    \ext{F_1} (\vb{x}) \join \ext{F_2} (\vb{x})$.
  \end{proof}

 \begin{thm}
   Let $m \in \N$, let $\ob{L}$ be a finite lattice.
   We suppose that we can compute the join and the meet of two
   elements in $\ob{L}$.
   Let $G$ be a finite subset of $\NW{m} \times \ob{L}$. Then
   there is an algorithm to determine whether $F_G$ satisfies
   \HC{7}.
 \end{thm}
 \begin{proof}
   For each triple $i,j,k$ with $\lambda_i \join \lambda_j = \lambda_k$,
   we have to check whether
          $F(\vb{a} + \vb{e}_k) = F(\vb{a} + \vb{e}_i) \join F(\vb{a} + \vb{e}_j)$ for all $\vb{a} \in \NW{m}$.
   For $\vb{x} \in \NW{m}$, let  $F_1 (\vb{x}) := F_G (\vb{x} + \vb{e_k})$, and let
   $F_{2} (\vb{x}) := F_G (\vb{x} + \vb{e}_i) \join F_G (\vb{x} + \vb{e}_j)$.
   Both $F_1$ and $F_2$ are antitone functions.

   We will now show that we can effectively evaluate
   $\ext{F_1}$ and $\ext{F_2}$ at each $\vb{a} \in \NWC{m}$.
   From Lemma~\ref{lem:shift}, we obtain
   $\ext{F_1} (\vb{a}) = \ext{F_G} (\vb{a} + \vb{e}_k)$.
   Therefore $\ext{F_G} (\vb{a} + \vb{e}_k)$ can be computed
   using  Lemma~\ref{lem:fga}\eqref{it:fg1}.
   Now let $F_3 (\vb{x}) := F_G (\vb{x} + \vb{e}_i)$
   and $F_4 (\vb{x}) := F_G (\vb{x} + \vb{e}_j)$.
   Then using Lemma~\ref{lem:shift}, we obtain $\ext{F_3} (\vb{a}) = \ext{F_G} (\vb{a} + \vb{e}_i)$
   and $\ext{F_4} (\vb{a}) = \ext{F_G} (\vb{a} + \vb{e}_j)$.
   Lemma~\ref{lem:join} yields that
   $\ext{F_2} (\vb{a}) = \ext{F_3} (\vb{a}) \join \ext{F_4} (\vb{a})$.
   Therefore
   $\ext{F_2} (\vb{a}) = \ext{F_G} (\vb{a} + \vb{e}_i) \join
   \ext{F_G} (\vb{a} + \vb{e}_j)$, and the evaluations of $\ext{F_G}$
   can again be accomplished using Lemma~\ref{lem:fga}\eqref{it:fg1}.

   By Theorem~\ref{thm:learn}, we can compute finite subsets
   $G_1$ and $G_2$ of $\NW{m} \times \ob{L}$
   such that $F_1 = F_{G_1}$ and $F_2 = F_{G_2}$.
   Now by two applications of Lemma~\ref{lem:compare}, we find
   out whether $F_{G_1} = F_{G_2}$.
 \end{proof}
 We notice that canonical representations of $F_1$ and $F_2$ can be computed
 more efficiently.
 In the proof of Theorem~\ref{thm:ccan}, we saw that we can easily find
 a canonical representation of a function $H$ if we
 know the minimal elements of
 $\{ \vb{a} \in \NW{m} \mid H (\vb{a}) \le \alpha \}$ for each $\alpha \in
 \ob{L}$. We will now see how to find these minimal elements for
 $F_1$ and $F_2$.
  Lemma~\ref{lem:ufg} provides us with a finite set
 $A$ such that
 $A{\uparrow} = \{ \vb{a} \in \NW{m} \mid F_G (\vb{a}) \le \alpha \}$.
 Then for each $\vb{b} \in \NW{m}$,
 we have $F_1 (\vb{b}) \le \alpha$ if and only if
 $\vb{b} + \vb{e}_k \in A{\uparrow}$. This is equivalent to
 $\vb{b} \in \{ \vb{a} - (\vb{a} \sqcap \vb{e}_k) \mid \vb{a} \in A \}{\uparrow}$.
 Thus we find the minimal elements of
 $\{ \vb{a} \in \NW{m} \mid F_1 (\vb{a}) \le \alpha\}$ as the
 minimal elements of
 $\{ \vb{a} - (\vb{a} \sqcap \vb{e}_k) \mid \vb{a} \in A \}$.
 Now we consider $F_2$. The minimal elements of
 $\{ \vb{a} \in \NW{m} \mid F_2 (\vb{a}) \le \alpha\}$ can
 be found using the following consideration:
 if we know that for two antitone functions $H_1$ and $H_2$,
 we have
 $\{ \vb{a} \in \NW{m} \mid H_1 (\vb{a}) \le \alpha \} = A_1{\uparrow}$
 and
 $\{ \vb{a} \in \NW{m} \mid H_2 (\vb{a}) \le \alpha \} = A_2{\uparrow}$,
 then
 $\{ \vb{a} \in \NW{m} \mid H_1 (\vb{a}) \join H_2 (\vb{a}) \le \alpha \} =
 (A_1{\uparrow}) \cap (A_2{\uparrow})$.
 In this way, we can find
 the minimal elements of $\{ \vb{a} \in \NW{m} \mid F_1 (\vb{a}) \le \alpha\}$ and
 $\{ \vb{a} \in \NW{m} \mid F_2 (\vb{a}) \le \alpha\}$ for each
 $\alpha \in \ob{L}$. Then we proceed as in the proof of Theorem~\ref{thm:ccan}
 to find canonical representations of $F_1$ and $F_2$.
In this way, we do not have to
 refer to Theorem~\ref{thm:learn}.

 In \cite{AM:SOCO, Mu:TLHC}, an operation sequence was called
 \emph{admissible} if it satisfies \HC1, \HC2, \HC3, \HC4, \HC7, and
 \HC8.  Now suppose that a sequence satisfies \HC3 and \HC4, and that
 $F_G$ is its encoding. We call $F_G$ \emph{admissible} if
 it satisfies \HC1, \HC2, \HC7, and \HC8.
 By combining the results from this section, we obtain:
 \begin{cor}
   Let $m \in \N$, and let  $\ob{L}$ be a finite lattice
   with $|\ob{L}| = m$
   in which we can compute joins and meets.
   Given a finite subset $G$ of $\NW{m} \times \ob{L}$,
   we can decide whether $F_G$ is admissible.
 \end{cor}
\section{Applications to commutator sequences} \label{sec:commseq}
In this section, we apply our results for representing
the sequence of commutator operations of a given algebra.
We will illustrate our result with two examples:
For $n \ge 2$,  we consider the algebras
$\ab{B}_n := \algop{\Z_4}{+, g_n}$ with
$g_n (x_1, \ldots, x_n) = 2 x_1 \cdots x_n$
and the algebra $\ab{B} := \algop{\Z_4}{+, (g_n)_{n \in \N}}$.
The polynomial functions of the algebras $\ab{B}_n$ provide an infinite ascending chain of clones
 on a four element set that was studied in \cite{Bu:OTNO}.
The congruence lattice of each of these algebras is the three
element chain $\{0,\alpha, 1\}$.
The commutator operations of $\ab{B}$
satisfy $[1,1]_{\ab{B}} = \alpha$, $[1, \alpha]_{\ab{B}} = 0$,
and $[1,\ldots,1]_{\ab{B}} = \alpha$ for every finite number of $1$'s.
The commutator operations of $\ab{C} := \ab{B}_7$ satisfy
$[1,\alpha]_{\ab{C}} = 0$, $[1,1]_{\ab{C}} = \alpha$,
and $[1,1,1,1,1,1,1]_{\ab{C}} = \alpha$,
$[1,1,1,1,1,1,1,1]_{\ab{C}} = 0$.

For an operation sequence $(f_n)_{n \in \N}$ on a lattice $\ob{L}$, we will
denote  $f_n (\alpha_1, \ldots, \alpha_n)$ by
$[\alpha_1, \ldots, \alpha_n]_f$. A \emph{commutator equality} over $\ob{L}$
is a formula $\varphi$ of the form
\[
  [\alpha_1, \ldots, \alpha_n] = \beta,
  \]
  where $n \in \N$, $\alpha_1, \ldots, \alpha_n, \beta \in \ob{L}$.
  We say that the operation sequence $(f_n)_{n \in \N}$ satisfies
  $\varphi$ if
  $[\alpha_1, \ldots, \alpha_n]_f = \beta$, and a universal
  algebra $\ab{A}$ satisfies $\varphi$ if its sequence of higher
  commutator operations
  $([.]_{\ab{A}}, [.,.]_{\ab{A}}, [.,.,.]_{\ab{A}}, \ldots)$ satisfies $\varphi$.

\begin{thm} \label{thm:largest}
  Let $\ab{A}$ be an algebra in a congruence modular variety
  with finitely many congruences.
  Then
  there is a finite set $\Phi$  of commutator equalities over its congruence lattice
  $\Con (\ab{A})$ such that the sequence of commutator operations of $\ab{A}$ is the largest
  sequence with \HC3 and \HC4 that satisfies all equalities in $\Phi$.
\end{thm}
\begin{proof}
  Let $m := | \Con (\ab{A}) |$, and let $\lambda_1, \ldots, \lambda_m$
  be the elements of $\Con (\ab{A})$.
  By Theorem~\ref{thm:tcadm}, the sequence of commutator operations of
  $\ab{A}$ satisfies \HC3 and \HC4. Therefore, as described in
  Section~\ref{sec:trans}, it can be encoded
  by an antitone function $F : \N_0^m \to \Con (\ab{A})$.
  By Theorem~\ref{thm:finrep1}, there is a finite subset of $G$
  of $F$ such that $F = F_G$.
  Now for every $((a_1, \ldots, a_m), \alpha) \in G$, we
  define the equality $\varphi_{(\vb{a}, \alpha)}$
  as
  \[
      [\, \underbrace{\lambda_1, \ldots, \lambda_1}_{a_1}, \ldots,
       \underbrace{\lambda_m, \ldots, \lambda_m}_{a_m}\,
      ]
      = \alpha.
   \]
   Then by Lemma~\ref{lem:largest},
   the sequence of commutator operations of $\ab{A}$ is the
   largest sequence satisfying \HC3 and \HC4 that satisfies
   all equalities in
   $\Phi_G := \{ \varphi_{(\vb{a}, \alpha)} \mid (\vb{a}, \alpha) \in G \}$.
\end{proof}
For example, the sequence of commutator operations of $\ab{B}$
is the largest sequence with \HC3 and \HC4 that
satisfies
\[
[1,1] = \alpha,
[1, \alpha] = 0,
[\alpha, \alpha] = 0,
[\alpha] = \alpha,
[0] = 0.
\]
Omitting those
 conditions that follow from
 (HC1) and (HC2), we obtain that
 the commutator sequence of $\ab{B}$
 is the largest operation sequence
  among those satisfying (HC1),(HC2),(HC3),(HC4)
with
$[1,1] = \alpha$ and $[1,\alpha] = 0$.

The sequence of commutator operations of $\ab{C}$
is the largest sequence with \HC3 and \HC4 that
satisfies
\[
  [1,1,1,1,1,1,1,1] = 0,
  [1,1] = \alpha,
  [1, \alpha]  = 0,
[\alpha, \alpha] = 0,
[\alpha] = \alpha,
[0] = 0.
\]
Omitting those
 conditions that follow from
 (HC1) and (HC2), we obtain that
 the commutator sequence of $\ab{B}$
 is the largest operation sequence
  among those satisfying (HC1),(HC2),(HC3),(HC4)
  with
  $[1,1] = \alpha$,
$[1,\alpha] = 0$,
and $[1,1,1,1,1,1,1,1] = 0$.

It follows from Proposition~\ref{pro:nonunique} that there is
no finite set $\Phi$ of commutator equalities such that the
sequence of commutator operations of $\ab{B}$ is uniquely
determined by $\Phi$ and the conditions \HC3 and \HC4.
In order to define such sequences uniquely, we define
\emph{extended commutator equalities}.
An \emph{extended commutator} over the complete lattice
$\ob{L}$ is an expression of the form
\[
    [S; \alpha_1, \ldots \alpha_n],
\]
where $n \in \N_0$, $S \subseteq \ob{L}$, and
$\alpha_1, \ldots, \alpha_n \in \ob{L}$.
For an operation sequence $(f_n)_{n \in \N}$ on $\ob{L}$ we define
its value by
\begin{equation} \label{eq:Salpha}
[S; \alpha_1, \ldots \alpha_n]_f :=
  \bigmeet \{ [\sigma_1, \ldots, \sigma_m, \alpha_1, \ldots,
    \alpha_n]_f  \, \mid \,
  m \in \N_0, (\sigma_1, \ldots, \sigma_m) \in S^m \};
\end{equation}
  additionally, for the empty sequence $\Lambda$  and $S = \emptyset$,
  $[\emptyset; \Lambda]$ is defined as the maximal element of $\ob{L}$.
  The following Lemma allows to omit those $\alpha_i$ that occur in $S$.
  \begin{lem}  \label{lem:omit}
    Let $(f_n)_{n \in \N}$ be an operation sequence on the complete
    lattice $\ob{L}$ with \HC3 and \HC4, let $S \subseteq \ob{L}$,  let $n \in \N$,
    and let $\alpha_1, \ldots, \alpha_n \in \ob{L}$. If
    $\alpha_1 \in S$, then
    $[S; \alpha_1, \alpha_2,  \ldots, \alpha_n]_f =[S; \alpha_2,  \ldots, \alpha_n]_f$.
  \end{lem}
  \begin{proof}
  Let
  $[\sigma_1, \ldots, \sigma_m, \alpha_1, \ldots,  \alpha_n]_f$ be one
  of the expressions in~\eqref{eq:Salpha}
  whose meet is $[S; \alpha_1, \alpha_2,  \ldots, \alpha_n]_f$.
  Then setting $\sigma_{m+1} := \alpha_1$, we see that the same expression
  appears as one of the expressions whose meet is
  $[S; \alpha_2,  \ldots, \alpha_n]_f$. This proves
    $[S; \alpha_1, \alpha_2,  \ldots, \alpha_n]_f \ge
     [S; \alpha_2,  \ldots, \alpha_n]_f$.
  For proving $\le$, let  $[\sigma_1, \ldots, \sigma_m, \alpha_2, \ldots,  \alpha_n]_f$ be one
  of the expressions in~\eqref{eq:Salpha}
  whose meet is $[S; \alpha_2,  \ldots, \alpha_n]_f$. Then since
  $(f_n)_{n \in \N}$ satisfies \HC3, we have
  $[\sigma_1, \ldots, \sigma_m, \alpha_2, \ldots,  \alpha_n]_f
  \ge
     [\alpha_1, \sigma_1, \ldots, \sigma_m, \alpha_2, \ldots,  \alpha_n]_f$,
     which by \HC4 is equal to
     $[\sigma_1, \ldots, \sigma_m, \alpha_1, \alpha_2, \ldots,  \alpha_n]_f$,
     which occurs as one of those expressions whose meet is
     $[S; \alpha_1, \alpha_2,  \ldots, \alpha_n]_f$.
     This proves   $[S; \alpha_2,  \ldots, \alpha_n]_f \ge
     [S; \alpha_1, \alpha_2,  \ldots, \alpha_n]_f$.
  \end{proof}
  We also observe that for a sequence $(f_n)_{n \in \N}$ with \HC4,
  we have $[S; \alpha_1, \ldots, \alpha_n] =
  [S; \alpha_{\pi(1)}, \ldots, \alpha_{\pi(n)}]$ for each
  permutation $\pi$ on $\ul{n}$.

  An \emph{extended commutator equality} over $\ob{L}$ is an expression
  $\varphi$ of the
form $[S; \alpha_1, \ldots, \alpha_n] = \beta$ with
$S \subseteq \ob{L}$, $n \in \N_0$, and
$\alpha_1, \ldots, \alpha_n, \beta \in \ob{L}$, and $(f_n)_{n \in \N}$
satisfies this equality if $[S; \alpha_1, \ldots, \alpha_n]_f = \beta$.
An   algebra $\ab{A}$ satisfies $\varphi$ if its sequence of higher
  commutator operations %
  $([.]_{\ab{A}}, [.,.]_{\ab{A}}, [.,.,.]_{\ab{A}}, \ldots)$ satisfies $\varphi$.
  We note that the  commutator equality $[\alpha_1, \ldots, \alpha_n] = \beta$ is
  equivalent to the extended equality $[\emptyset; \alpha_1, \ldots, \alpha_n]
  = \beta$.
  The following Lemma tells how $[S; \alpha_1, \ldots, \alpha_n]$
  is linked to the encoding of $(f_n)_{n \in \N}$.
  \begin{lem} \label{lem:encodingOfExtension}
    Let $\ob{L}$ be a finite lattice, let $(f_n)_{n \in \N}$ be
    an operation sequence on $\ob{L}$ with
    \HC3 and \HC4, and let $F : \NW{m} \to \ob{L}$ be its encoding.
    Let $S$ be a nonempty subset of $\ob{L}$, and let $\alpha_1, \ldots, \alpha_n \in \ob{L}$. For each $j \in \ul{m}$, let
    $a_j = \infty$ if $\lambda_j \in S$,
    and $a_j = |\{ k \in \ul{n} \,:\, \alpha_k = \lambda_j \}|$ if
    $\lambda_j \not\in S$.
    Then
    $[S; \alpha_1, \ldots, \alpha_n] =
    \ext{F} (a_1, \ldots, a_m)$.
  \end{lem}
  \begin{proof}
    Let $k \in \N$ be such that
    $|\ob{L} \setminus S| = k$, and let $t_1 < \cdots < t_k$ be
    such that $\ob{L} \setminus S = \{ \lambda_{t_1}, \ldots, \lambda_{t_k} \}$.
   Then
    \begin{multline*}
      \ext{F} (\vb{a})
      = \bigmeet \{ F (\vb{x}) \mid \vb{x} \in \NW{m}, \vb{x} \le \vb{a} \} \\
      = \bigmeet \{ F (\vb{x}) \mid \vb{x} \in \NW{m}, (x_{t_1}, \ldots, x_{t_k}) = (a_{t_1}, \ldots, a_{t_k}) \},
    \end{multline*}
    where the last equality follows from the fact that $F$ is antitone.
    From the fact that $(f_n)_{n \in \N}$ satisfies \HC4, we have
    \begin{multline*}
      \bigmeet \{ F (\vb{x}) \mid \vb{x} \in \NW{m}, (x_{t_1}, \ldots, x_{t_n}) = (a_{t_1}, \ldots, a_{t_n}) \} \\
     = \bigmeet \{ [\sigma_1, \ldots, \sigma_k, \underbrace{\lambda_{t_1}, \ldots, \lambda_{t_1}}_{a_{t_1}}, \ldots,
       \underbrace{\lambda_{t_k}, \ldots, \lambda_{t_k}}_{a_{t_k}}]_f \mid k \in \N_0, \sigma_1, \ldots, \sigma_k \in S \} \\
    = [S; \underbrace{\lambda_{t_1}, \ldots, \lambda_{t_1}}_{a_{t_1}}, \ldots,
      \underbrace{\lambda_{t_k}, \ldots, \lambda_{t_k}}_{a_{t_k}}]_f.
    \end{multline*}
    Now we use Lemma~\ref{lem:omit}
    in order
    to insert those $\alpha_i$ with
    $\alpha_i \in S$
and the symmetry of
the extended commutator coming from \HC4
and obtain
    $[S; \underbrace{\lambda_{t_1}, \ldots, \lambda_{t_1}}_{a_{t_1}}, \ldots,
      \underbrace{\lambda_{t_k}, \ldots, \lambda_{t_k}}_{a_{t_k}}]_f
    =[S; \alpha_1, \ldots, \alpha_n]_f$.
  \end{proof}
  For an $m$-element lattice $\ob{L} = \{\lambda_1, \ldots, \lambda_m\}$,
  $\vb{a} \in \NWC{m}$ and $\alpha \in \ob{L}$, we define
  an extended commutator equality $\varphi_{(\vb{a}, \alpha)}$ as follows:
  we set
  \begin{equation} \label{eq:ST}
     S :=  \{ \lambda_i \mid i \in \ul{m}, a_i = \infty \}, \,\,
     T  := \{ j \in \ul{m} \mid a_j \in \N_0 \},
  \end{equation}
  we let $k := |T|$  and we assume  that $t_1 < \cdots < t_k$ are the elements
  of $T$.
  Then
      $\varphi_{(\vb{a}, \alpha)}$ is defined as the extended equality
  \[
    [S;
        \underbrace{\lambda_{t_1}, \ldots, \lambda_{t_1}}_{a_{t_1}}, \ldots,
        \underbrace{\lambda_{t_k}, \ldots, \lambda_{t_k}}_{a_{t_k}}
    ] = \alpha.
    \]
  \begin{lem} \label{lem:transeeq}
    Let $\ob{L}$ be a finite lattice with $m$ elements,
    let $(f_n)_{n \in \N}$ an operation sequence on $\ob{L}$ that
    satisfies \HC3 and \HC4, let $\vb{a} \in \NWC{m}$, and
    let $\alpha \in \ob{L}$.
    Then $(f_n)_{n \in \N}$ satisfies $\varphi_{(\vb{a}, \alpha)}$
    if and only if the extension of its encoding $F$ as an antitone function
    satisfies
    $\ext{F} (\vb{a}) = \alpha$.
  \end{lem}
  \begin{proof}
    Let $S, T$ be as in~\eqref{eq:ST}. Then the value of
    the left hand side of $\varphi_{(\vb{a}, \alpha)}$ is
    $[S;
        \underbrace{\lambda_{t_1}, \ldots, \lambda_{t_1}}_{a_{t_1}}, \ldots,
        \underbrace{\lambda_{t_k}, \ldots, \lambda_{t_k}}_{a_{t_k}} ]_f$,
    which by Lemma~\ref{lem:encodingOfExtension} is equal to
    $\ext{F} (a_1, \ldots, a_m)$.
 \end{proof}
\begin{thm} \label{thm:unique}
  Let $\ab{A}$ be an algebra in a congruence modular variety with
  finitely many congruences.
  Then
  there is a finite set $\Phi$  of extended commutator equalities over its congruence lattice
  $\Con (\ab{A})$ such that its sequence of commutator operations is the
  unique operation
  sequence on $\Con (\ab{A})$ with \HC3 and \HC4 that satisfies all extended equalities in $\Phi$.
\end{thm}
\begin{proof}
    Let $m := | \Con (\ab{A}) |$, and let $\lambda_1, \ldots, \lambda_m$
  be the elements of $\ob{L} := \Con (\ab{A})$.
  By Theorem~\ref{thm:tcadm}, the sequence of commutator operations of
  $\ab{A}$ satisfies \HC3 and \HC4. Therefore, as described in
  Section~\ref{sec:trans}, it can be encoded
  by an antitone function $F : \NW{m} \to \Con (\ab{A})$.
  Let $\ext{F}$ be its extension to $\NWC{m}$ as defined in~\eqref{eq:Fhat}.
  Then Theorem~\ref{thm:hatF} yields a finite subset $H$ of $\ext{F}$ such that
  $\ext{F}$ is the unique antitone function from $\NWC{m}$ to $\ob{L}$
  that contains $H$ as a subset.
  We define
  $$\Phi := \{ \varphi_{(\vb{a}, \alpha)} \mid (\vb{a}, \alpha) \in H \}.$$
  By Lemma~\ref{lem:transeeq},
  the sequence of commutator operations of $\ab{A}$ satisfies $\Phi$. To
  show uniqueness,
  let $(f_n)_{n \in \N}$ be an operation sequence
  with \HC3 and \HC4 that satisfies $\Phi$, and let $F_1$ be its
  encoding.
  From Lemma~\ref{lem:transeeq}, we obtain that
   $H \subseteq \ext{F_1}$. Since $\ext{F}$ is the only
  antitone function with $H \subseteq \ext{F}$, we
  have $\ext{F} = \ext{F_1}$. Therefore, $(f_n)_{n \in \N}$ has
  $\ext{F}$ as the extension of its encoding, and thus
  $(f_n)_{n \in \N}$ is equal to
  the sequence of higher commutator operations of $\ab{A}$.
\end{proof}

 For the algebra $\ab{B}$, we
 obtain
 \begin{multline*}
  H  = \{ ((0,0,2), \alpha), ((0,1,1), 0), ((0,2,0), 0),
  ((1,0,0), 0), ((0,1,0), \alpha),  \\
  ((0,0,1), 1),
  ((0,0,\infty), \alpha),
 ((\infty, \infty, \infty),0) \}
  \end{multline*}
  as a complete representation of the encoding of its sequence of commutator operations.
  These translate into the extended commutator equalities
  \begin{multline*}
    [1,1] = \alpha, [1, \alpha] = 0, [\alpha, \alpha] = 0, [0] = 0,
    [\alpha] = \alpha, \\
      [1] = 1,
      [\{1\}; \Lambda ] = \alpha,
      [ \{0, \alpha, 1\} ; \Lambda] = 0,
  \end{multline*}
   where $\Lambda$ denotes the empty sequence.
 By omitting those conditions that follow from
 (HC1) and (HC2), we obtain that
 the commutator sequence of $\ab{B}$
 is the unique operation sequence
  among those satisfying (HC1),(HC2),(HC3),(HC4)
with
$[1,1] = \alpha$,  $[1,\alpha] = 0$, $[\alpha ] = \alpha$,
$[1] = 1$, $[\{1\}; \Lambda] = \alpha$; omitting the conditions
of the form $[\sigma] = \sigma$, we obtain that
 the commutator sequence of $\ab{B}$
 is the unique operation sequence
 among those satisfying (HC1),(HC2),(HC3),(HC4) and
 $[\sigma] = \sigma$ for all $\sigma \in \{0,\alpha,1\}$ that
 satisfies
$[1,1] = \alpha$,  $[1,\alpha] = 0$ and $[\{1\}; \Lambda] = \alpha$.

Similarly, we obtain that the commutator sequence of
$\ab{B}_7$ the unique one  among those with (HC1), (HC2), (HC3), (HC4),
and $[\sigma] = \sigma$ for all $\sigma \in \{0, \alpha, 1\}$ that satisfies
 $[1,\alpha] = 0$, $[1,1] = \alpha$, $[1,1,1,1,1,1,1] = \alpha$,
 and $[1,1,1,1,1,1,1,1] = 0$.

 Now suppose that we are given an algebra $\ab{A}$ in a congruence
 modular variety with finitely many congruences, and that we can compute its higher commutators
 \[
   [\alpha_1, \ldots, \alpha_k]_{\ab{A}}
 \]
 and
 \[
   [S; \alpha_1, \ldots, \alpha_k]_{\ab{A}}
    :=   \bigmeet \{ [\sigma_1, \ldots, \sigma_m, \alpha_1, \ldots,
    \alpha_n]_{\ab{A}}  \, \mid \,
   m \in \N, \sigma_1, \ldots, \sigma_m \in S \}
   \]
   for all $k \in \N$, $n \in \N_0$ and for all nonempty subsets $S$
   of $\Con (\ab{A})$.
   Then Theorem~\ref{thm:learn} yields that we can learn the complete
   sequence of higher commutator operations of $\ab{A}$ from
   finitely many evaluations of the form
   $[S; \alpha_1, \ldots, \alpha_k]_{\ab{A}}$:
   \begin{thm} \label{thm:learn2}
     Let $\ab{A}$ be an algebra in a congruence modular variety
     with finitely many congruences.
     Then, evaluating only finitely  many expressions of
     the form
     \[
       [S; \alpha_1, \ldots, \alpha_k]_{\ab{A}}
       \]
     with $S \subseteq \Con (\ab{A})$ and
     $\alpha_1, \ldots, \alpha_n \in \Con (\ab{A})$,
     we can compute a finite set $\Phi$ of extended commutator equalities
     such that
     the sequence of higher commutator operations on $\ab{A}$ is
     the unique operation sequence on $\Con (\ab{A})$ with
     \HC3 and \HC4 that satisfies $\Phi$.
   \end{thm}
   \begin{proof}
     We use Theorem~\ref{thm:learn} to obtain a representation
     of the encoding of the sequence of commutator operations of
     $\ab{A}$, and  Theorem~\ref{thm:ccomplete} to obtain a complete
     representation. Then translating the complete representation
     into extended commutator equalities, we obtain the required set $\Phi$.
   \end{proof}
   We say that an algebra with finitely many congruences
   has a \emph{computable extended commutator
     sequence} if we can compute with the elements of its congruence lattice
   as described in Section~\ref{sec:compcan}, and if  there is
   furthermore an algorithm that, given
   $S \subseteq \Con (\ab{A})$ and $\alpha_1, \ldots, \alpha_n \in
   \Con (\ab{A})$,
   computes $[S; \alpha_1, \ldots, \alpha_n]_{\ab{A}}$.
   For such an algebra, we are also able to check
   whether it satisfies a given extended commutator equality
   $[S; \alpha_1, \ldots, \alpha_n] = \beta$.
   \begin{cor}
     There is an algorithm that,
     given two algebras
     $\ab{A}_1, \ab{A}_2$ in a congruence modular variety
     on the same universe with the same finite congruence lattice
     and
     computable extended commutator sequences, decides
     whether $\ab{A}_1$ and $\ab{A}_2$ have the same sequence
     of higher commutator operations.
   \end{cor}
   \begin{proof}
     We use Theorem~\ref{thm:learn2} to compute a set
     $\Phi_2$ of extended commutator equalities
     which uniquely determine the higher commutator operations
     of $\ab{A}_2$. We return ``true'' if
     sequence of  higher commutator operations of $\ab{A}_1$ satisfies
     $\Phi_2$.
   \end{proof}

   \bibliography{axiomComm43}
\end{document}